\documentclass[onefignum,onetabnum]{siamart171218}

\usepackage{lipsum}
\usepackage{amsfonts}
\usepackage{amsmath}
\usepackage{graphicx}
\usepackage{epstopdf}
\usepackage{algorithmic}
\usepackage{amssymb}
\usepackage{enumitem}
\usepackage{bbm}
\usepackage[ruled,vlined,algo2e]{algorithm2e}
\usepackage{xcolor}
\usepackage{caption}
\ifpdf
  \DeclareGraphicsExtensions{.eps,.pdf,.png,.jpg}
\else
  \DeclareGraphicsExtensions{.eps}
\fi

\newsiamremark{remark}{Remark}
\newsiamremark{hypothesis}{Hypothesis}
\crefname{hypothesis}{Hypothesis}{Hypotheses}
\newsiamthm{claim}{Claim}
\newsiamremark{assumption}{Assumption}
\headers{A Globally Convergent LP/SOCP Algorithm for SDP}{B. Roig-Solvas and M. Sznaier}

\title{A Globally Convergent LP and SOCP-based algorithm for Semidefinite Programming}

\author{Biel Roig-Solvas\thanks{Northeastern University, Boston MA
  }
\and Mario Sznaier\footnotemark[1]}

\usepackage{amsopn}

\makeatletter
\newcommand*{\addFileDependency}[1]{
  \typeout{(#1)}
  \@addtofilelist{#1}
  \IfFileExists{#1}{}{\typeout{No file #1.}}
}
\makeatother

\newcommand*{\myexternaldocument}[1]{%
    \externaldocument{#1}%
    \addFileDependency{#1.tex}%
    \addFileDependency{#1.aux}%
}

\ifpdf
\hypersetup{
  pdftitle={A Globally Convergent LP/SOCP Algorithm for SDP},
  pdfauthor={B. Roig-Solvas and M. Sznaier}
}
\fi

\myexternaldocument{ex_supplement}

\begin{document}

\maketitle

\begin{abstract}
  Semidefinite programs (SDP) are one of the most versatile frameworks in numerical optimization, serving as generalizations of many conic programs and as relaxations of NP-hard combinatorial problems. Their main drawback is their computational and memory complexity, which sets a practical limit to the size of problems solvable by off-the-shelf SDP solvers. To circumvent this fact, many algorithms have been proposed to exploit the structure of particular problems and increase the scalability of SDPs for those problem instances. Progress has been less steep, however, for general-case SDPs. In this paper, motivated by earlier results by Ahmadi and Hall, we show that a general SDP can be solved to $\epsilon$-optimality, in polynomial time, by performing a sequence of less computationally demanding Linear or Second Order Cone programs. In addition, we provide a bound on the number of iterations required to achieve $\epsilon$-optimality. These results are illustrated using random SDPs and well-known problems from the SDPLib dataset.

\end{abstract}



\section{Introduction}\label{Section:Introduction}
Semidefinite programming (SDP) is one of the most versatile frameworks in the field of convex optimization and encompasses a wide range of problems in science and engineering \cite{Vander1996,SDPHB2000}. In addition to serving as a generalization of conic problems like linear programming (LP), quadratic (QP) and second-order cone programming (SOCP), they provide convex relaxations for a large set of non-convex problems including combinatorial, rank-constrained and polynomial optimization problems \cite{SDPHB2011}. The general standard form SDP is of the form:\\
\begin{equation}
    \begin{split}
        X_{PSD}^* = \operatorname*{argmin}_X &\quad \operatorname{Tr}\left(C^T\,X\right)\\
        \text{s.t.}&\quad \operatorname{Tr}\left(A_i^T\,X\right) = b_i\quad i= 1,\dots,M\\
        & \quad X \in  \mathcal{S}^+_N
    \end{split}
    \label{1:eq:SDP}
\end{equation}
where $\mathcal{S}^+_N$ denotes the cone of $N\times N$ positive semi definite (PSD) matrices.  The flexibility of SDPs  comes, however, at a computational cost. For instance, interior point solvers scale  as $\mathcal{O}\left(M N^3 + M^2 N^2\right)$ \cite{alizadeh1998}, limiting applicability to relatively small problems.  This issue becomes more pressing when SDPs are used as convex relaxations of non-convex problems, e.g. polynomial optimization, which increases the dimensionality of the problems and worsens its scaling. Several approaches have been proposed throughout the last decades to deal with this computational barrier, exploiting  sparsity \cite{Fujisawa1997,Fukuda2001,Nakata2003,Andersen2010,Vander2015} or  alternative implementations relying on first-order methods \cite{Zheng2017} or non-convex programming \cite{Burer2003,Burer2005}. An alternative approach seeks to obtain lower complexity relaxations by replacing the semi-definite constraints with less computationally demanding  linear or second order cone constraints \cite{Ahmadi2015,Ahmadi2019} , leading to an algorithm that alternates between Cholesky decompositions and linear (LP) or second order cone (SOCP) programs.  However, while successful in many scenarios, there is no guarantee that it will converge to the solution of the original SDP. Indeed, there are examples where the approach  fails to attain the optimum.
 Thus, the issue of whether a general SDP can be solved via LP or SOCP was left open.

This paper gives an affirmative answer to this question. Our main result shows that a generic SDP can be solved to $\epsilon$-optimality in polynomial time by performing a sequence of LPs and SOCPs, combined with Cholesky factorizations. Further, we provide an upper bound on the number of iterations that depends only on the problem data. The key observation is that optimality can be guaranteed by periodically returning to the central path of the original SDP, a problem that can also be solved using LP/SOCPs. Thus, the proposed algorithm iterates between decreasing the objective function and visiting the central path until eventually converging to an $\epsilon$-optimizer of the SDP, measured by the SDP's duality gap. To the best of our knowledge, this is the first globally convergent SDP algorithm based on DD and SDD programs.

The paper is organized as follows. Section \ref{Section:PrevWork} provides background material on DD and SDD programs and Interior Point Methods (IPMs). The proposed algorithm and its proof of convergence are presented in Sections  \ref{Section:PropAlgo} and \ref{Section:Convergence}. We present numerical results in Section \ref{Section:NumericalP} and discuss complexity considerations and algorithm extensions in Section \ref{Section:Discussion}. Finally, directions for future work  are given in Section \ref{Section:Conclusion}.

\section{Previous Work}\label{Section:PrevWork}
In this section  we recall some theoretical properties of IPMs that will be used in the analysis of the proposed algorithm and briefly summarize the main work our approach rests on,  the diagonally-dominant (DD) and scaled diagonally-dominant (SDD) relaxations of SDPs  \cite{Ahmadi2015,Ahmadi2019}.

\subsection{Interior Point Methods for SDPs and the Central Path}

Interior point methods (IPM) are arguably the most common algorithms for general purpose convex conic programs, including LP, SOCP and SDPs. First formulated  in 1984  \cite{karmarkar1984}, IPMs have become widely adopted due to their guaranteed polynomial runtime  \cite{nesterov1994}. These methods  handle conic constraints by adding to the cost a  penalty function (a ``barrier")  that tends to infinity when approaching the boundary of the feasible set from inside the set.
To prevent numerical instability, IPMs solve a sequence of optimization problems in which the barrier is weighted by a factor $1/t$, where $t$ is
increased until $\epsilon$-optimality is reached. In the case of SDPs, the most widely used barrier function is the negative log-determinant, which leads to problems of the form:
\begin{equation}
    \begin{split}
        X_t = \operatorname*{minimize}_X &\quad \operatorname{Tr}\left(C^T\,X\right) - \frac{1}{t}\log\left(|X|\right)\\
        \text{s.t.}&\quad \operatorname{Tr}\left(A_i^T\,X\right) = b_i\quad i= 1,\dots,M
    \end{split}
    \label{4:eq:BarrierPr}
\end{equation}

The curve in 
$\mathcal{S}^+$ defined by the optimizers $X_t$ of \eqref{4:eq:BarrierPr} as a function of $t>0$ is called the Central Path of the problem. As  $t\to\infty$, $X_t$ converges to $X_{PSD}^*$, the optimizer of  \eqref {1:eq:SDP}, i.e. $X_{PSD}^*$. Moreover, due to duality theory, the elements of the central path satisfy the following inequality:
\begin{equation}
 \operatorname{Tr}\left(C^T\,X_t \right)\geq \operatorname{Tr}\left(C^T\,X^*_{PSD}\right)\geq \operatorname{Tr}\left(C^T\,X_t\right)-N/t   
 \label{2:eq:dualbound}
\end{equation}
which provides an optimality bound at any point in the path. Further, it can be shown that, given a strictly increasing sequence $t_k$, the corresponding cost sequence is strictly decreasing, e.g. $\operatorname{Tr}\left(C^T\,X_t\right)<\operatorname{Tr}\left(C^T\,X_{t-1}\right)$.

\subsection{Optimization of self-concordant functions via IPMs}\label{sec:selconcordant}

The proof of convergence of the proposed algorithm rests on the  properties of the central path for self-concordant functions, e.g. those satisfying $|f^{'''}(x) \leq 2f^{''}(x)^\frac{3}{2}$ \cite{Boyd2004}. Consider the minimization of $f(x)$  through a Newton method \cite{Boyd2004}. Denote by $x,x_+$ and $x^*$  the current iterate, the iterate after taking a Newton step from $x$ and the global minimizer of $f(x)$, respectively,  and by $\lambda(x)$ the Newton decrement of $f(x)$ evaluated at $x$.  If the line-search constants are chosen such that $\alpha\in(0,0.5)$ and $\beta\in(0,1)$ and the variable $\eta = \frac{1-2\alpha}{4}$, then \cite{Boyd2004}:
\begin{equation}
     f(x_+) \leq f(x) - \alpha \beta \frac{ \lambda^2(x)}{1+ \lambda(x)}
    \label{4:eq:sc1}
\end{equation}
\begin{equation}
    \lambda(x) \leq \eta \quad \implies \quad f(x_+) \leq f(x) - \alpha \lambda^2(x)
    \label{4:eq:sc2}
\end{equation}
\begin{equation}
    \lambda(x) \leq 0.68 \quad \implies f(x) \geq f(x^*) \geq f(x) - \lambda^2(x)
    \label{4:eq:sc3}
\end{equation}

\subsection{DD and SDD relaxations of Semidefinite Programming}\label{2:Subsec:DD}
In \cite{Ahmadi2019}, the authors proposed a relaxation for general semidefinite programs based on replacing  the positive semidefinite constraints 
by lower complexity ones involving diagonally-dominant and scaled diagonally-dominant matrices, defined below:
\begin{definition}
A symmetric matrix $X$ is diagonally-dominant (DD) if\\
\begin{equation*}
    X(i,i)\;\geq \;\sum_{j\neq i}\;|X(i,j)|\quad \forall i
    \label{2:def:dd}
\end{equation*}
\end{definition}
\begin{definition}
A symmetric matrix $X$ is scaled diagonally-dominant (SDD) if there exist a positive diagonal matrix $D$ and a DD matrix $Y$ such that $X = D Y D$.
\label{2:def:sdd}
\end{definition}
From Gershgorin circle theorem it follows that DD and SDD matrices are positive semidefinite. Further, the following inclusion holds   $DD_N \subset SDD_N \subset PSD_N$, where  $DD_N,SDD_N$ and $PSD_N$ denote the cones of $N\times N$ DD, SDD and PSD matrices. Thus, relaxations of the SDP
 \eqref {1:eq:SDP} can be obtained by  replacing the constraint $X\in \mathcal{S}_N^+$ with the stronger ones $X \in DD_N$ or $X\in SDD_N$. 
The following resuls, adapted from  \cite{Ahmadi2015,Boman2005} provides an alternative characterization of DD and SDD matrices that was used in  \cite{Ahmadi2015} to show that these relaxations lead to lower complexity LPs or SOCPs.
 Define the mapping $\Psi_{i,j}$ from $2\times2$ matrices to $N\times N$ matrices:
\[
   \Psi_{i,j}(M) = \bar{M}
   \quad\text{where}  
  \left \{ \begin{tabular}{c} $\bar{M}(\{i,j\},\{i,j\}) = M$  \\
  $0$ \text{ otherwise}.
  \end{tabular} \right.
\]
i.e. the $\{i,j\}$ sub-matrix of $\bar{M}$ is $M$, and the rest of entries of $\bar{M}$ are 0 \footnote{Whenever necessary, if the first argument of $\Psi(\cdot,\cdot)$ contains a set of subindices $i,j$, we will omit the second argument, i.e.  $\Psi\left(M_{i,j},\{i,j\}\right) = \Psi\left(M_{i,j}\right)$.}.  $\Psi_{i,j}(.)$ allows for characterizing the set of $DD$ and $SDD$ matrices in terms of ``exploded" 2$\times$2 matrices as follows:

\begin{lemma}[\cite{Ahmadi2015,Boman2005}]  \label{lem:SDD}
\[ Y\in DD_N \iff Y = \sum_{i,j}^N \Psi_{i,j} (M_{i,j}), \quad M_{i,j} \in DD_2 \]
Similarly,
\[Y\in SDD_N \iff Y = \sum_{i,j}^N \Psi_{i,j}(M_{i,j}), \quad M_{i,j} \succeq 0 \]
\end{lemma}
From this result it follows that enforcing the constraint $X \in DD_N(SDD)_N$ indeed reduces to a set of linear (second order cone) constraints.

\subsection{Iterative Basis Update}\label{2:Subsec:IBU}
Replacing the PSD constraint in \eqref {1:eq:SDP} with the stronger one $X \in DD_N$ or $X\in SDD_N$ leads to a computationally cheaper optimization. However, the solution to these relaxed problems can be far from the true optimum. To address this,  \cite{Ahmadi2015} proposed an iterative algorithm, based on alternating between  solving a sequence of DD/SDD problems and performing Cholesky factorizations. Briefly, the idea is to solve at step $k$ a problem of the form
\begin{equation}
    \begin{split}
     &  X_k^*\left(U_{k-1}\right) =  \operatorname*{argmin}_{X,Y}\; \operatorname{Tr}\left(C^TX\right)\\
      &  \text{s.t.}\; \operatorname{Tr}\left(A_i^TX\right) = b_i\quad i=1,\dots,M\\
        &U_{k-1}^T Y U_{k-1} = X,\quad Y\in DD_N\;/\;SDD_N
    \end{split}
    \label{2:eq:IBUDDSDD}
\end{equation}
where $U_{k-1}$ is a Cholesky factor of the previous solution, e.g. $X_{k-1}=U_{k-1}^TU_{k-1}$. Since  $I \in DD_N (SDD_N)$,  the previous iterate $ X_{k-1}^*$, is always a feasible solution of \eqref{2:eq:IBUDDSDD}. Hence the algorithm generates a sequence of solutions $X_i^*$, with non-increasing cost.
This sequence, however, is not guaranteed to converge to the optimizer of the  SDP \eqref{1:eq:SDP} and in numerical tests tends to converge to strictly suboptimal values for all medium to large size problems ($N>>10$). The proposed algorithm avoids becoming trapped in these suboptimal accumulation points by periodically ``centering'' the iterates by projecting onto the central path of \eqref{1:eq:SDP}.

\section{Proposed Algorithm}\label{Section:PropAlgo}

In this section we present the proposed algorithm for solving the SDP \eqref{1:eq:SDP} to $\epsilon$-optimality. The algorithm is split in two phases. The first phase,  the \emph{decrease phase}, consists of solving a sequence of DD/SDD programs, exactly as in \cite{Ahmadi2015}. As noted above, this sequence  tends to stagnate on a suboptimal objective cost as the iterates approach the boundary of the PSD cone and their conditioning worsens. To prevent this, a second phase of the algorithm starts after the decrease phase that consists of a series of steps designed to improve the iterates' conditioning. We call these steps \emph{centering} steps, as they guide the iterates towards the center path of the SDP by solving a sequence of analytic centerings on the DD/SDD set. These centering steps constitute the \emph{centering phase} of the algorithm.

\begin{minipage}{0.95\linewidth}
\makebox[\linewidth]{
 \includegraphics[width = 6cm]{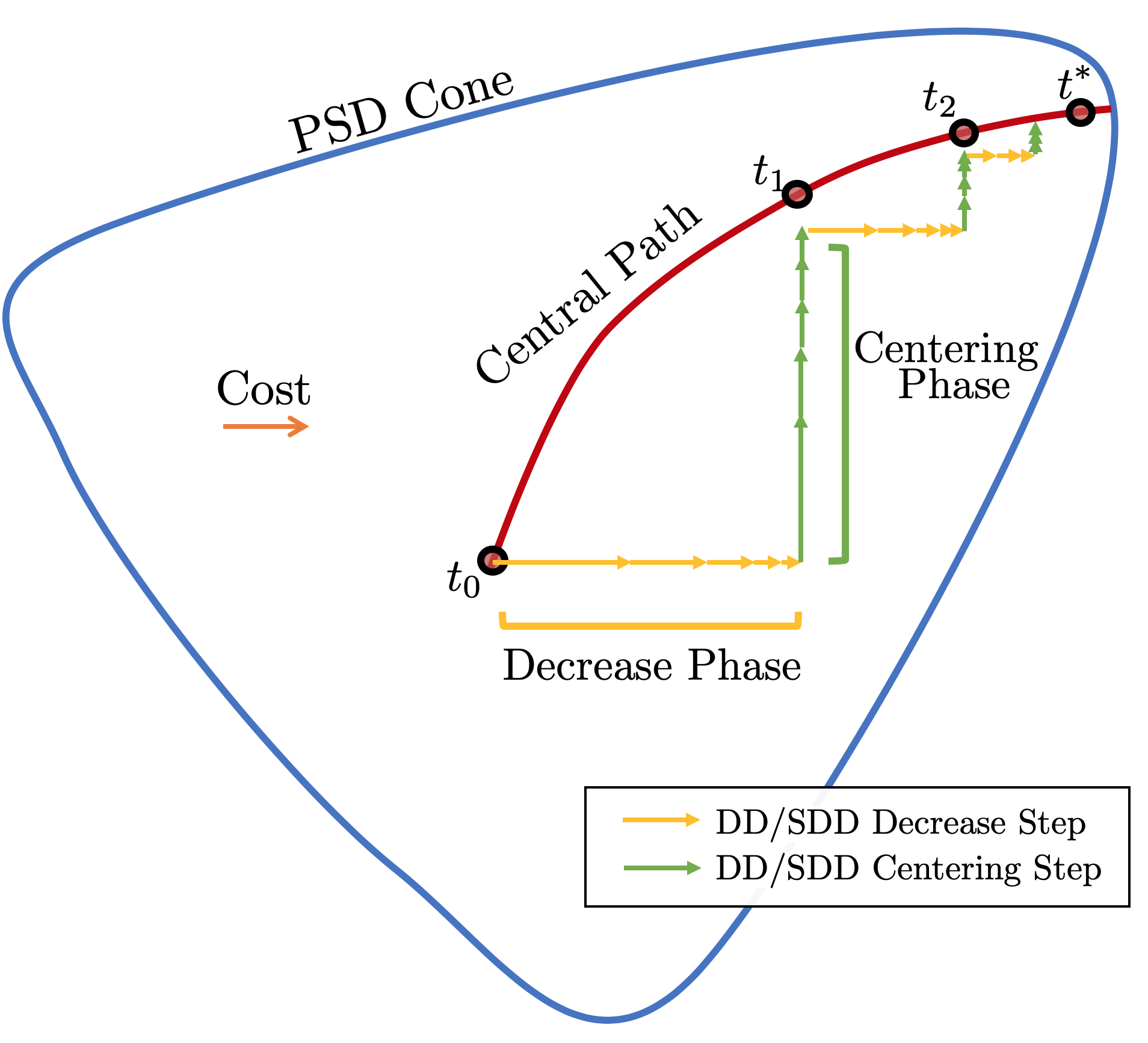}}
\captionof{figure}{Sketch of the proposed algorithm. The algorithm alternates between  a \emph{decrease phase}, where the cost is decreased, and a
\emph{centering phase} that brings the iterate to a point $\epsilon_c$-close to the central path. After  $\kappa$ decrease and centering phases, the iterate reaches a point $\epsilon_c$-close to the central path with parameter $t_{\kappa} \geq t^*$, guaranteeing $\epsilon_g$-convergence to the global optimizer of the SDP.} \label{3:fig:algorithmsketch}
\end{minipage}

Figure \ref{3:fig:algorithmsketch} illustrates  the proposed algorithm. In the \emph{decrease phase}, a sequence of problems of the form \eqref{2:eq:IBUDDSDD} are solved,  decreasing the cost. After a given number  $s_d$ of decrease steps\footnote{A study on the impact the choice $s_d$ has on the algorithm performance is carried out in Section \ref{Section:NumericalP}. The theoretical guarantees developed in Section \ref{Section:Convergence} hold for any value of $s_d\geq1$, and so we assume $s_d=1$ unless otherwise stated.}, the \emph{centering phase} starts and a sequence of  analytic centering problems of the form \eqref{3:eq:centering} are solved :
\begin{equation}
    \begin{split}
       X_{l} \left(U_{l-1}\right) =  \operatorname*{argmin}_{X,Y} &\;-\phi\left(Y\right)\\
        \text{s.t.}\; &\operatorname{Tr}\left(A_i^TX\right) = b_i\quad i=1,\dots,M\\
        &\operatorname{Tr}\left(C^TX\right) = \operatorname{Tr}\left(C^TX_{l-1}\right)\\
        &U_{l-1}^T Y U_{l-1} = X,\quad Y\in DD_N\;/\;SDD_N
    \end{split}
    \label{3:eq:centering}
\end{equation}
where the function $-\phi(Y)$ is the logarithmic barrier of the DD/SDD sets and $U_{l-1}$ is the Cholesky factor of $X_{l-1}$. The sequence of centering steps converges to a point $\epsilon_c$ close to the central path of the SDP, whose  optimality gap can be found from \eqref{2:eq:dualbound}. At this point, a  new \emph{decrease phase} starts and the algorithm keeps alternating between  \emph{decrease} and \emph{centering} phases, as outlined in Algorithm 1,  until it converges to an  $\epsilon_g$-optimal  solution of the original problem  \eqref{1:eq:SDP}.  In the next section  we will prove that the algorithm converges  in polynomial time, and provide a bound on the number of iterations as a function of $\epsilon_g$ and the problem data.

\begin{algorithm}[H]
\SetAlgoLined
\KwResult{$\epsilon_g$-optimal $X^*$ }
 Initialize optimality gap as $g_{SDP} = \infty$; Initialize $X_0$; $k=0$;\\
 \While{$g_{SDP} > \epsilon_g$ }{
  \% Start Decrease Phase\\
  \For{Number of decrease steps $s_d\geq1$}{
  $k = k+1$; Compute Cholesky Factors of $X_{k-1}$;\\
  Solve for $X_{k}(U_{k-1})$ as in \eqref{2:eq:IBUDDSDD};\;
  }
  Update optimality gap $g_{SDP}$;\\
   \% Start Centering Phase\\
  Initialize centering gap $g_C = \infty$; $l=0$, $X_{k,l} = X_k$;\\ 
  \While{$g_C > \epsilon_C$}{
  $l = l+1$; Compute Cholesky Factors of $X_{k,l-1}$;\\
   Solve for $X_{k,l}(U_{k,l-1})$ as in  \eqref{3:eq:centering};\\
   Update centering gap $g_C$;\;
   }
   Update optimality gap $g_{SDP}$;\;
   $X_k = X_{k,l}$;
 }
 $X^* = X_k$;\;
 \caption{Globally Convergent DD/SDD Algorithm}\label{3:Algo:Principal}
\end{algorithm}

\section{Convergence}\label{Section:Convergence}

In this section we prove global convergence of Algorithm \ref{3:Algo:Principal}. The proof  relies on (i) establishing a polynomial upper bound on the number of instances that problem \eqref{3:eq:centering} needs to be solved to achieve  $\epsilon_c$-optimality;  and  (ii) a  proof of convergence of the combined Decrease and Centering Phases to an $\epsilon_g$-optimizer of the SDP \eqref{1:eq:SDP} also in polynomial time. We start with the assumptions that will hold for both proofs:
\begin{assumption}
The data matrices $C$ and $A_i$ all satisfy that:
\begin{equation*}
\begin{split}
    \operatorname{Tr}\left(A_i^TA_j\right) = \operatorname{Tr}\left(A_i^TC\right) = 0 \quad &\forall 1\leq i,j\leq M \\
    ||C||_F = ||A_i||_F = 1 \quad &\forall i = 1,\dots,M
\end{split}
\end{equation*}
\label{4:ass:1}
\end{assumption}
\begin{assumption}
\eqref{1:eq:SDP} admits a feasible  $X\succ 0$ (Slater's condition).
\label{4:ass:2}
\end{assumption}
\begin{assumption}
The optimizer of \eqref{1:eq:SDP} satisfies $\operatorname{Tr}(C^TX^*_{PSD}) > -\infty$.
\label{4:ass:3}
\end{assumption}

Assumption \ref{4:ass:1} can be made to hold trivially for any SDP,  by orthogonalizing the  matrices $A_i$ and projecting out the component of $C$ spanned by these matrices. Assumption \ref{4:ass:2} is required to guarantee strong duality (see \cite{permenter2018partial,pataki2013strong} and references therein for face reduction methods to deal with problems with no strictly feasible solution) and Assumption \ref{4:ass:3} guarantees that the optimal cost function is finite.
\subsection{Convergence of the Centering Phase} \label{Subsec:Centering}
The goal of this proof is to show that the sequence ${X}_l$ in \eqref{3:eq:centering} converges to the optimizer of the PSD analytic centering:
\begin{equation}
    \begin{split}
        X\left( U_{k-1}\right) = \operatorname*{argmin}_{X,Y} \; -h(Y), \quad
        \text{s.t.}\; &\operatorname{Tr}\left(A_i^TX\right) = b_i\quad i=1,\dots,M\\
        &\operatorname{Tr}\left(C^TX\right) = \operatorname{Tr}\left(C^TX_{k-1}\right)\\
        &U_{k-1}^T Y U_{k-1} = X,\quad Y \succeq 0
    \end{split}
    \label{4:eq:PSDcent}
\end{equation}
where the objective function is defined as:
\begin{equation}
    h(Y) = \left(N-1\right) \log\left(|Y|\right) - N \left(N-1\right)\log\left(N-1\right)
\end{equation}
i.e. a scaled and shifted variant of the common log-determinant barrier for the PSD cone. Motivated by Lemma \ref{lem:SDD}, we will consider the following
 logarithmic barriers:
\begin{equation}
\begin{split}
    \phi_{DD}(\mathcal{M}) &= \frac{1}{2}\sum_{i,j>i}^N \log\left(M_{i,j}(1,1)^2- M_{i,j}(1,2)^2\right) \\&+ \log\left(M_{i,j}(2,2)^2- M_{i,j}(1,2)^2\right)\\
    \phi_{SDD}(\mathcal{M}) &= \sum_{i,j>i}^N \log\left(M_{i,j}(1,1)\,M_{i,j}(2,2) - M_{i,j}(1,2)^2\right) 
\end{split}
\label{4:eq:logbarrs}
\end{equation}

The next Lemma establishes key properties of these barrier functions.
\begin{lemma}\label{4:lemma:selfconcor}
$-\phi_{DD}(\mathcal{M})$ and $-\phi_{SDD}(\mathcal{M})$  have the following properties:
\begin{enumerate}
\item[(a)] Self-concordance with respect to the entries of $\mathcal{M}$.
\item[(b)] If $Y = \Psi(\mathcal{M}) \in DD_N$, then $-\phi_{DD}(\mathcal{M})\geq-\phi_{SDD}(\mathcal{M})$.
\item[(c)]
If $Y = \Psi(\mathcal{M}) \in SDD_N$, then $-\phi_{SDD}(\mathcal{M})\geq-h(Y)$.
\end{enumerate}
\end{lemma}
\begin{proof}
See  Section \ref{App:PC:selfconcor}
\end{proof}
\begin{corollary}
If $Y = \Psi(\mathcal{M}) \in DD_N$, then $-\phi_{DD}(\mathcal{M})\geq-\phi_{SDD}(\mathcal{M})\geq-h(Y)$. Moreover, if $M_{i,j} = \frac{1}{N-1} I_{2\times2}$ for all $i<j$, then $Y = \Psi(\mathcal{M}) =I$ and  $-\phi_{DD}(\mathcal{M})=-\phi_{SDD}(\mathcal{M})=-h(Y) = N \left(N-1\right)\log\left(N-1\right)$
\label{4:coro:Id}
\end{corollary}
\begin{proof}
The first statement follows from properties (b) and (c). The second follows from evaluating $\phi_{DD}(\mathcal{M})$ and $\phi_{SDD}(\mathcal{M})$
at $M_{i,j} = \frac{1}{N-1} I_{2\times2}$ and $-h(Y)$ at $Y=I$. 
\end{proof}

The convergence proof proceeds by comparing the evolution of a Newton method applied to problems \eqref{3:eq:centering} and \eqref{4:eq:PSDcent}.  Assume one is solving problem \eqref{3:eq:centering} using a Newton method where the variable $Y$ is parametrized as $Y = \Psi(\mathcal{M})$. Denote by $\mathcal{M}_0$ the set where all elements of the set are $M_{i,j} = \frac{1}{N-1} I_{2\times 2}$, i.e. $Y_0 = \Psi(\mathcal{M}_0) = I$, which by construction is always a feasible solution of \eqref{3:eq:centering}. Denoting the evaluation of the Newton decrement for $-\phi(\mathcal{M} = \mathcal{M}_0)$\footnote{As the proof applies to both DD and SDD cases, we drop the subscripts for notation clarity.} by  $\lambda_{\phi}(\mathcal{M}_0)$, then properties \eqref{4:eq:sc1}--\eqref{4:eq:sc3} lead to the following two lemmas: 
\begin{lemma}
Let $U_{l-1}$ to be the Cholesky factor of $X_{l-1}$. Then if the Newton decrement of problem \eqref{3:eq:centering} satisfies $\lambda_{\phi}(\mathcal{M}_0) > \frac{\eta}{\sqrt{N-1}}$, its optimizer $X_l(U_{l-1})$ satisfies:
\begin{equation*}
    -(N-1)\log\left(|X_l|\right) \leq -(N-1)\log\left(|X_{l-1}|\right) - \xi
\end{equation*}
where $\xi$ is a positive constant of the form:
\begin{equation*}
    \xi = \frac{\alpha\beta}{\sqrt{N-1}}\,\frac{\eta^2}{\sqrt{N-1}+\eta}
\end{equation*}
\label{4:lemma:SCdamped}
\end{lemma}
\begin{proof}
See proof in Section \ref{App:PC:SCdamped}.
\end{proof}
\begin{lemma}
If the Newton decrement of \eqref{3:eq:centering} satisfies $\lambda_{\phi}(\mathcal{M}_0) \leq \frac{\eta}{\sqrt{N-1}}$, and the centering optimality gap between $X_{l-1}$ and the optimizer $X_*$ of \eqref{4:eq:PSDcent} is given by:
\begin{equation*}
    g_{l-1}= (N-1)\,\left(\log\left(|X_*|\right) -\log\left(|X_{l-1}|\right)\right) 
\end{equation*}
then the centering optimality gap $g_l$ between $X_l$ and $X_*$ is upper-bounded as:
\begin{equation*}
   \frac{N-1-\alpha}{N-1} \min(\eta^2,g_{l-1}) \geq g_l \geq (N-1)\,\left(\log\left(|X_*|\right) -\log\left(|X_l|\right)\right) 
\end{equation*}
\label{4:lemma:SCquad}
\end{lemma}
\begin{proof}
See proof in Section \ref{App:PC:SCquad}.
\end{proof}
Lemmas \ref{4:lemma:SCdamped} and \ref{4:lemma:SCquad} provide the foundation for the proof of polynomial complexity of the Centering Phase, as shown next:
\begin{theorem}
The Centering Phase described in Algorithm \ref{3:Algo:Principal} converges to $\epsilon_C$-optimality to the optimizer of \eqref{4:eq:PSDcent} as:
\begin{equation*}
   \epsilon_C \geq  (N-1)\,\left(\log\left(|X_*|\right) -\log\left(|X_{L}|\right)\right)
\end{equation*}in at most $L$ iterations, where $L$ is given by:
\begin{equation*}
    L = \lceil\frac{(N-1)\,\left(\log\left(|X_*|\right) -\log\left(|X_{0}|\right)\right)-\epsilon_C}{\xi}\rceil + \lceil \frac{\log\left(\epsilon_C\right)-\log\left(\eta^2\right)}{\log\left(N-1-\alpha\right)-\log\left(N-1\right)}\rceil
\end{equation*}
ans $X_0$ is the starting point of the Centering Phase.
\end{theorem}
\begin{proof}
The proof follows easily from the results in Lemmas \ref{4:lemma:SCdamped} and \ref{4:lemma:SCquad}. At each iteration the centering optimality gap is reduced, either by a fixed amount $\xi$ if $\lambda_{\phi}(\mathcal{M}_0) > \frac{\eta}{\sqrt{N-1}}$, as shown in the first Lemma, or by a multiplicative factor $\frac{N-1-\alpha}{N-1}$ if $\lambda_{\phi}(\mathcal{M}_0) \leq \frac{\eta}{\sqrt{N-1}}$, as given by the second Lemma. Bringing the centering gap below $\epsilon_C$ requires at most $L_1$ iterations for the fixed decrease, with $L_1 = \lceil\left((N-1)\,\left(\log\left(|X_*|\right) -\log\left(|X_{0}|\right)\right)-\epsilon_C\right)/\xi\rceil $, or $L_2$ iterations for the relative decrease, where $ L_2 = \lceil \left(\log\left(\epsilon_C\right)-\log\left(\eta^2\right)\right)/\left(\log\left(N-1-\alpha\right)-\log\left(N-1\right)\right)\rceil$, leading to a total running time of at most $L1+L2$ iterations:
\begin{equation}
\begin{split}
    L = L1 + L2 = \lceil\frac{(N-1)\,\left(\log\left(|X_*|\right) -\log\left(|X_{0}|\right)\right)-\epsilon_C}{\xi}\rceil \\ + \lceil \frac{\log\left(\epsilon_C\right)-\log\left(\eta^2\right)}
    {\log\left(N-1-\alpha\right)-\log\left(N-1\right)}\rceil
    \end{split}
    \end{equation}
\end{proof}

\subsection{Convergence of the Proposed Algorithm Under Perfect Centering} \label{Subsec:Conv:Dec}
In this section we develop the proof of $\epsilon_g$-convergence of Algorithm \ref{3:Algo:Principal}.
 The main idea  is to show that alternating between decreasing and centering phases leads to a sequence of solutions $X_k$ which are identical to the ones  obtained using an interior point algorithm to  solve \eqref{4:eq:BarrierPr} for a specific sequence $t_k$ that satisfies
$t_k > \chi t_{k-1}$, where $\chi >1$ is a constant that depends on the problem data. It follows that a desired value $t^*$ (corresponding to a given optimality gap) can be found in at most $\kappa  = \lceil \frac{\log\left(N/\epsilon^*\right) - \log\left(t_0 \right)}{\log\left(\chi\right)}\rceil$ iterations.  For simplicity,  we assume here that $\epsilon_C = 0$, i.e. exact convergence of the \emph{centering} phase.  Then, in Section \ref{Subsec:Conv:Inexact} we adapt these results for the practical case of $\epsilon_C >0$. 
  The proof is divided into the following steps. 
\begin{enumerate}
\item \textbf{Step 1:} Show that after each combined decrease and centering steps, the cost $\operatorname{Tr}\left(C^T {X}_k \right)$ decreases at least by an amount  $\frac{\sqrt{\Phi}}{t_{k-1}\sqrt{1+\Theta}}$, where $\Phi$ and $\Theta$ depend only on the problem data. 
 \item \textbf{Step 2:} Use the result above, combined with the strictly decreasing property of  the solution to  \eqref{4:eq:BarrierPr} to establish that $t_k > \chi t_{k-1}$ for some constant $\chi > 0$.  Thus we can reach any desired $t^*$ in at most $\kappa = \log(t^*/t_0)/\log(\chi)$ iterations. Using \eqref{2:eq:dualbound}, this value of $t
^*$ translates to an optimality gap $\epsilon^* = N/t^*$.
\end{enumerate}

 \noindent \textbf{Step 1:} We start the proof by recasting problem \eqref{4:eq:PSDcent} into its simplest form:
\begin{equation}
    \begin{split}
       \tilde{ X}_c = \operatorname*{argmin}_{X\succeq 0} &\quad -\log\left(|X|\right)\\
        \text{s.t.}&\quad \operatorname{Tr}\left(A_i^T\,X\right) = b_i\quad i= 1,\dots,M\\
        & \quad \operatorname{Tr}\left(C^T\,X\right) = c
    \end{split}
    \label{4:eq:PSDcent2}
\end{equation}
and showing that the optimizer of \eqref{4:eq:PSDcent2} is also an optimizer of \eqref{4:eq:BarrierPr} for a specific value of $t$:
\begin{lemma}
Given $\tilde{X}_c$ the solution of problem \eqref{4:eq:PSDcent2}, the solution $X_t$ of \eqref{4:eq:BarrierPr} coincides with $\tilde{X}_c$ when $t = |\tau|$, where $\tau$ is the dual variable of \eqref{4:eq:PSDcent2} associated to the linear constraint $\operatorname{Tr}\left(C^T\,X\right) = c$ evaluated at $X = \tilde{X}_c$.
\label{4:lemma:CentEquiv}
\end{lemma}

\begin{proof} Follows  from the KKT stationarity conditions of problems \eqref{4:eq:PSDcent2} and \eqref{4:eq:BarrierPr}.
\end{proof}
Next, we use the result above to obtain a  (potentially conservative)
 bound on the cost decrease after a \emph{decrease} step is taken.
 
 \begin{lemma}\label{4:lemma:decrease}
 Let $X_{k}$ denote the  optimizer of the Centering Phase  at stage $k$.  Consider a positive constant $\Phi$ that satisfies $\Phi \leq \{2/(N+1),1\}$ for the DD and SDD cases, respectively. Define $\hat{Y}_k \doteq I - \sqrt{\Phi}\frac{Q}{||Q||_F}$, where $Q \doteq U_{k-1}^{-T} C U_{k-1}^{-1}$, and $\hat{X}_k \doteq U_{k-1}^T \hat{Y}_k U_{k-1}$. Then (a) $\hat{X}_k$ is a feasible solution to  \eqref{2:eq:IBUDDSDD}, and (b)
  $\operatorname{Tr}\left(C^T X_{k}\right)
     \leq \operatorname{Tr}\left(C^T \hat{X}_{k}\right) \leq \operatorname{Tr}\left(C^T X_{k-1}\right) - \frac{\sqrt{\Phi}}{t_{k-1}\sqrt{1+\Theta}}$, where $\Theta$ is a finite positive constant that depends only on the problem data. 
\end{lemma}
\begin{proof} Given in Section \ref{App:PD:decrease} \end{proof}

\noindent \textbf{Step 2:} From the strictly decreasing property of the objective of \eqref{4:eq:BarrierPr} with respect to $t$ (Lemma \ref{App:PD:alb}), it follows that the values of $t$ such that the corresponding solution to  \eqref{4:eq:BarrierPr} are  $X_k$ and $X_{k-1}$ satisfy $t_k > t_{k-1}$. However, in order to establish finite time convergence we need to prove that  $t_k \geq  \chi t_{k-1}$,  for some  $\chi>0$. In turn, this requires, given $X_k$, finding the corresponding value $t_k$. Since this is a non-trivial problem, we will find a function $g(.)$ such that $g(t_k) \leq   \operatorname{Tr}\left(C^T X_k\right)$ and use it as a proxy, to find some $\tilde{t_k} \leq t_k$. The desired results will be established by showing that
$t_k \geq \tilde{t_k} \geq  \chi t_{k-1}$.  
\begin{lemma} Given $t_{k-1}$, let $X_{t_{k-1}}$ denote the solution to  \eqref{4:eq:BarrierPr} corresponding to $t=t_{k-1}$. Then
the  function $g_{k-1}(t) = \operatorname{Tr}\left(C^T X_{t_{k-1}}\right) - N/t_{k-1} + N/t$ is a lower bound of $\operatorname{Tr}\left(C^T X_{t}\right)$ for any $t\geq t_{k-1}$.
\label{4:lemma:gt}
\end{lemma}
\begin{proof} See Section \ref{App:PD:gt} \end{proof}
Choosing $\tilde{t}_k$ such that  $g_{k-1}(\tilde{t}_k )=  \operatorname{Tr}\left(C^T \hat{X}_k\right )$ we have that
\begin{equation}
\begin{split}
   \operatorname{Tr}\left(C^T {X}_{t_{k-1}}\right ) - \frac{N}{t_{k-1}} + \frac{N}{\tilde{t}_k}&  =  \operatorname{Tr}\left(C^T \hat{X}_k\right )  \leq  \operatorname{Tr}\left(C^T {X}_{t_{k-1}}\right )- \frac{\sqrt{\Phi}}{t_{k-1}\sqrt{1+\Theta}} \\
    &\implies \tilde{t}_k \geq t_{k-1} \frac{N \sqrt{1+\Theta}}{N \sqrt{1+\Theta} - \sqrt{\Phi}} \doteq t_{k-1} \; \chi
\end{split}
\label{4:eq:gchi}
\end{equation}
where $\chi >1$. Since
 $\operatorname{Tr}\left(C^T X_{t_k}\right)
     \leq \operatorname{Tr}\left(C^T \hat{X}_{k}\right) = g_{k-1}(\tilde{t}_k)  \leq \operatorname{Tr}\left(C^T X_{\tilde{t}_k}\right )$
and the cost is a decreasing function of $t$ it follows that $t_k \geq\hat{t}_k \geq \tilde{t}_k \geq  \chi t_{k-1}$,  leading to the main result of this section:
\begin{theorem} In the ideal case where the centering problem  \eqref {4:eq:PSDcent} can be solved exactly,
Algorithm \ref{3:Algo:Principal} converges to $\epsilon_g$-optimality in at most $\kappa$ iterations of Decrease and Centering Phases, where $\kappa$ is given by:
\begin{equation}
    \kappa  = \lceil \frac{\log\left(N/\epsilon^*\right) - \log\left(t_0 \right)}{\log\left(\chi\right)}\rceil\quad \text{with}\quad \chi = \frac{N \sqrt{1+\Theta}}{N \sqrt{1+\Theta} - \sqrt{\Phi}}>1
\label{4:eq:kappa}    
\end{equation}
\label{4:thm:dec_conv}
\end{theorem}
\begin{proof}
The proof follows  from equation \eqref{4:eq:gchi}. Starting from $t_0$, after each Centering and Decrease phase we increase the value of $t$ by a multiplicative factor greater than $\chi$, leading to $t_k > \chi\; t_{k-1}$. It follows that the number of iterations needed to reach $t^* = N/\epsilon^*$ is given by $t^* \leq \chi^\kappa \;t_0$ with $\kappa$ given by \eqref{4:eq:kappa}.
\end{proof}

 \subsection{Convergence of the Proposed Algorithm for $\epsilon_C>0$} \label{Subsec:Conv:Inexact}
In the previous section we have shown  convergence of the algorithm under the assumption that each Centering Phase ends with an iterate $X_{k,L}$ on the central path, i.e. $X_{k-1,L} = X_{t_{k-1}}$. In practice, however, we will not be able to reach the central path to exact precision and instead, through self-concordance, we can only guarantee that we are within a bound $\epsilon_C $ to the optimizer, in the form of:
\begin{equation}
       \epsilon_C \geq  (N-1)\,\left(\log\left(|X_{t_{k-1}}|\right) -\log\left(|X_{k-1,L}|\right)\right) \geq 0
       \label{4:eq:epsbound}
\end{equation}
Next we show that the convergence results in Theorem  \ref{4:thm:dec_conv}  hold even if the \emph{centering} phase provides an $\epsilon_c$ suboptimal solution $X_L$ to \eqref{4:eq:PSDcent}, provided that $\epsilon_c < \bar{\epsilon}_c$, where the constant $\bar{\epsilon}_c$ depends only on the problem data.  
The intuition behind the proof, illustrated in 
 Figure \ref{fig:fig3}, is that if $X_L$ and $X_{k}$\footnote{In the sequel, for notational simplicity, we refer to $X_{k-1,L}$ as $X_L$.}
 are close enough, then the solution $\hat{X}$ can be shown to satisfy the norm constraint $||U_L^{-T}\,\hat{X}\, U_L^{-1}-I||_F^2 \leq \Phi^\prime$, where $U_L$ are the Cholesky factors of $X_L$ and $\Phi
^\prime$ is a constant such that $\Phi < \Phi^\prime \leq \{2/(N+1),1\}$. Thus the $\hat{X}$ used in the proof of step 1 in section  \ref{Subsec:Conv:Dec}
 is within the feasible set of the decrease step \eqref{2:eq:IBUDDSDD} taken from $X_L$, rather than  $X_{k}$ (green circle in Fig \ref{fig:fig3}). 
To formalize these ideas we will proceed along the following steps:
\begin{enumerate}
\item \textbf{Step 3:} Find an upper bound of
$||U_L^{-T}\,\hat{X}\, U_L^{-1}-I||_F^2$ in terms of $\epsilon_C$, the centering optimality gap.
\item \textbf{Step 4:} Use this bound to show that, if  $\epsilon_c$ is within a specific range, then
  $\hat{X} =U_{k-1}^T \hat{Y} U_{k-1}$  is still a feasible solution to \eqref{2:eq:IBUDDSDD}, when this step is started  from $X_L=U_L^TU_L$ rather than  $X_{t_{k-1}}$. Hence, the decreasing step yields a cost no larger than $\operatorname{Tr}\left(C^T \hat{X} \right)$, proving convergence of the algorithm for the suboptimal case $\epsilon_C>0$ (see Fig. \ref{fig:fig3}).
\end{enumerate}

\begin{minipage}{0.9\linewidth}
\makebox[\linewidth]{
 \includegraphics[width = 6cm]{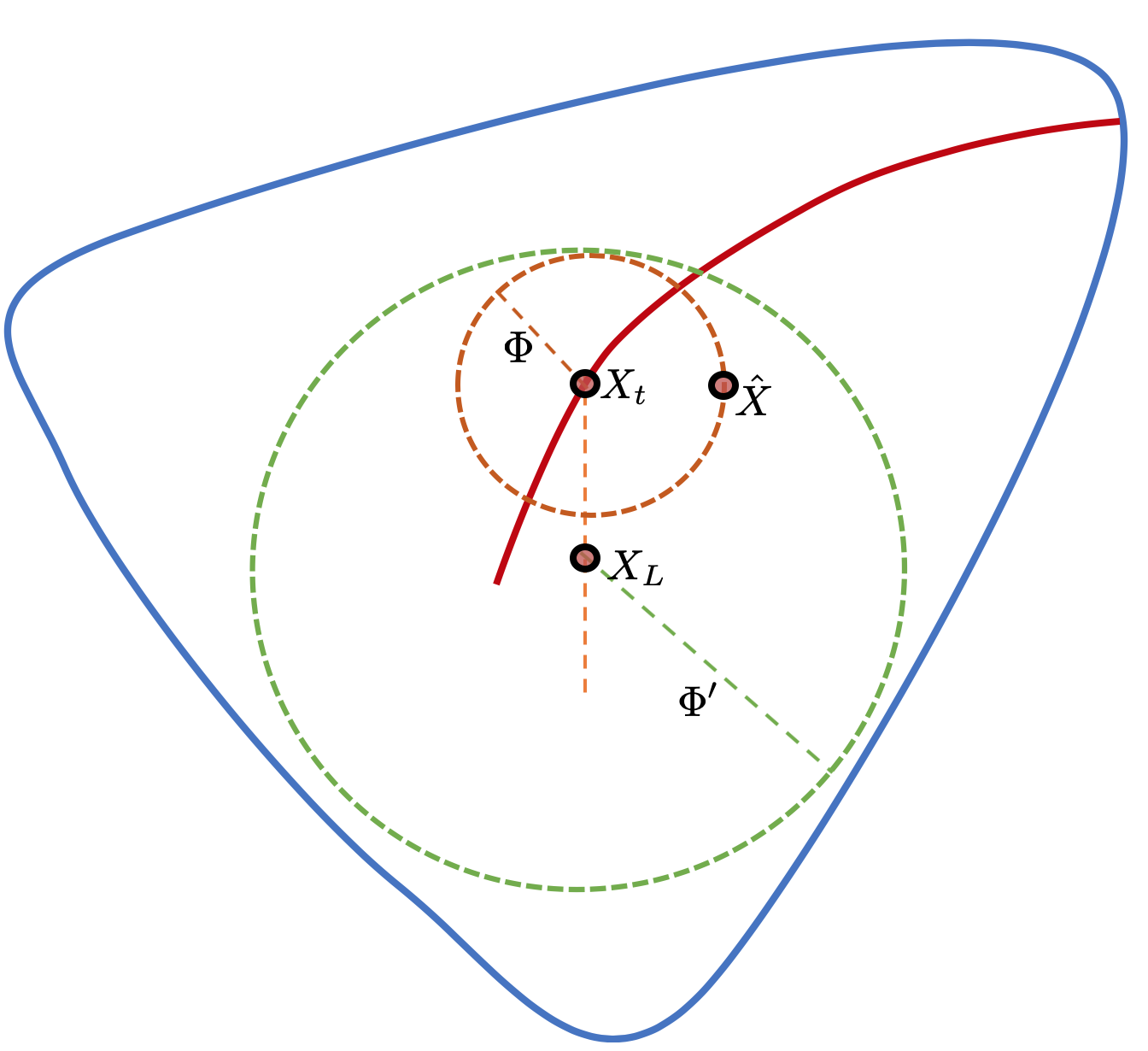}}
\captionof{figure}{Outline of the proof of convergence when the \emph{centering phase} stops at a $\epsilon_C$-optimal point $X_L$. Circles in orange and green indicate volumes that lie within a decrease step of $X_t$ and $X_L$ , respectively. If $X_t$ and $X_L$  are close enough, measured by $\epsilon_C$, then the orange volume is contained within the green one and $\hat{X}$ lies within a decrease step from $X_L$, extending the guarantees of Theorem 2 to the suboptimal $\epsilon_C>0$ case.} \label{fig:fig3}
\end{minipage}

\noindent \textbf{Step 3:} The next result provides the required bound:
\begin{lemma}\label{4:lemma:ineq}
 The following inequality holds:
\begin{equation}
    ||U_L^{-T}\left(\hat{X}-X_L\right) U_L^{-1}||_F^2 \leq \Phi^\prime
    \label{4:eq:ellipgoal}
\end{equation}
\end{lemma}
\begin{proof} Start by considering the chain of inequalities:
\begin{equation}
    \begin{split}
    &  ||U_L^{-T}\left(\hat{X}-X_L\right) U_L^{-1}||_F =||U_L^{-T}\left(\hat{X}-X_t+X_t-X_L\right) U_L^{-1}||_F\\
      &\leq||U_L^{-T}\left(\hat{X}-X_t\right) U_L^{-1}||_F +||U_L^{-T}\left(X_t-X_L\right) U_L^{-1}||_F\\
      &=||U_L^{-T}U_t^{T}U_t^{-T}\left(\hat{X}-X_t\right)U_t^{-1}U_t U_L^{-1}||_F+||U_L^{-T}\left(X_t-X_L\right) U_L^{-1}||_F\\
      &\leq||U_L^{-T} X_t U_L^{-1} ||_2\;||U_t^{-T}\left(\hat{X}-X_t\right)U_t^{-1}||_F +\sqrt{N}||U_L^{-T}\left(X_t-X_L\right) U_L^{-1}||_2\\
      &\leq||U_L^{-T} X_t U_L^{-1} ||_2\;\sqrt{\Phi} +\sqrt{N}||U_L^{-T}\left(X_t-X_L\right) U_L^{-1}||_2
    \end{split}
    \label{4:eq:inexactineqs}
\end{equation}
where the first inequality is due to the triangle inequality, the second is due to the inequalities between the Frobenius and the spectral norm given by $||A^T B A ||_F^2\leq ||A^TA||_2^2 ||B||_F^2$ and $||A||_F^2 \leq N \,||A||_2^2$ and the third one due to $||U_t^{-T}\,\hat{X}\, U_t^{-1}-I||_F^2 \leq \Phi$. The last two  spectral norms can be bounded using the results in \cite{Roig2020}
by:
\[ \begin{aligned}
    ||U_L^{-T} X_t U_L^{-1} ||_2 & \leq \frac{1}{1-\sqrt[3]{1-e^{-\epsilon_C}}},\\
    ||U_L^{-T}\left(X_t-X_L\right) U_L^{-1}||_2& \leq \frac{\sqrt[3]{1-e^{-\epsilon_C}}}{ 1-\sqrt[3]{1-e^{-\epsilon_C}} }
    \end{aligned}
    \]
  leading to:
  \begin{equation}
    ||U_L^{-T}\left(\hat{X}-X_L\right) U_L^{-1}||_F \leq \frac{\sqrt{\Phi}}{1-\sqrt[3]{1-e^{-\epsilon_C}}} +  \frac{\sqrt{N} \, \sqrt[3]{1-e^{-\epsilon_C}}}{ 1-\sqrt[3]{1-e^{-\epsilon_C}} }
\end{equation} 
The right hand side is an increasing function of $\epsilon_C$: it reduces to $\sqrt{\Phi}$ for $\epsilon_C = 0$ and both denominators tend to 0 as $\epsilon_C$ grows. To certify that equation \eqref{4:eq:ellipgoal} holds, it suffices to upper bound the right hand side by $\sqrt{\Phi^\prime}$. To  find the largest $\epsilon_C$  that makes this bound tight, define $z = \sqrt[3]{1-e^{-\epsilon_C}}$, leading to:
\begin{equation}
    \frac{\sqrt{\Phi}}{1-z} +  \frac{\sqrt{N}\,z}{ 1-z } = \sqrt{\Phi^\prime}\implies z = \frac{ \sqrt{\Phi^\prime}-\sqrt{\Phi}}{\sqrt{N}+\sqrt{\Phi^\prime}}
    \label{4:eq:zPhi}
\end{equation}
where $z>0$ due to the assumption $\Phi^\prime > \Phi$. \end{proof}
\noindent \textbf{Step 4}:
From the result above and Lemma  \ref{4:lemma:sphere}, it follows that if \eqref{4:eq:zPhi} holds, then
$Y_L = U_L^{-T}\hat{X}U_L^{-1} \in DD/SDD$ and hence $\hat{X}$ is indeed a feasible solution of  \eqref{2:eq:IBUDDSDD}, when starting  from $X_L=U_L^TU_L$.
This observation leads to the following Theorem:
\begin{theorem}
Assume $X_L$ is an $\epsilon_C$-optimal solution to the Centering Phase. Assume further that $\Phi^\prime = \{2/(N+1),1\}$, for the DD and SDD case respectively, and that $\Phi<\Phi^\prime$. Then if:
\begin{equation}
    \epsilon_C \leq    \epsilon_c^* = -\log\left(1- \left( \frac{ \sqrt{\Phi^\prime}-\sqrt{\Phi}}{\sqrt{N}+\sqrt{\Phi^\prime}}\right)^3\right)
\end{equation}
the cost function after a decrease step taken from $X_L$ is at least as low as the cost function evaluated at $\hat{X}$.
\end{theorem}
\begin{proof}
 Follows from undoing the change of variables $z = \sqrt[3]{1-e^{-\epsilon_C}}$ in equation \eqref{4:eq:zPhi}. This choice of $\epsilon_C$ guarantees that $||U_L^{-T}\left(\hat{X}-X_L\right) U_L^{-1}||_F^2 \leq \Phi^\prime$ and thus $\hat{X}$ lies within the feasible space of a decrease step taken from $X_L$. Thus, the optimizer of the decrease step will have an objective cost at least as low as that of $\hat{X}$.
\end{proof}

  \subsection{Termination criteria} \label{Subsec:Conv:FindingT}
The estimation of the central path parameter $t$ is necessary for the termination of the algorithm,  as it provides a way to compute the duality gap $\epsilon = \frac{N}{t}$ and thus an estimate on how suboptimal the current iterate is.  In subsection \ref{Subsec:Conv:Dec} we have shown that the convergence of the proposed algorithm to any target $t^*$ is guaranteed, and complemented this result in subsection \ref{Subsec:Conv:Inexact} with the proof that this convergence still holds even when terminating the Centering phase at a suboptimal iterate.

However, while the parameter $t$ can be recovered from the dual variables of the centering problem at its optimum, as shown in Lemma \ref{4:lemma:CentEquiv}, the same does not hold if the centering problem is terminated at a suboptimal iterate $X_{k,L}$. Next we show that we can obtain computable bounds on the parameter $t$ for the suboptimal termination case and use those to determine the duality gap of the iterate $X_{k,L}$ using only variables available to the algorithm at that execution point. 

\begin{theorem}
Assume $X_{k,L}$ is an $\epsilon_C$-optimal solution to the Centering Phase. Then the parameter $t_k$ associated with the exact optimizer of that Centering Phase is bounded above and below by:
\begin{equation}\begin{aligned}
  t_k &\leq  -\tau +  \left(\left(1+\frac{1}{\alpha\beta}\right)\epsilon_c + \sqrt{2\epsilon_c}\right)||X_{k,L}^{-1}||_F \; \text{and} \\
  t_k &\geq \;-\tau -  \left(\left(1+\frac{1}{\alpha\beta}\right)\epsilon_c + \sqrt{2\epsilon_c}\right)||X_{k,L}^{-1}||_F
\end{aligned}
\end{equation}
where $\tau$ is the value of the dual variable associated to the constraint $\operatorname{Tr}(C^TX) = \operatorname{Tr}(C^TX_{l-1})$ in \eqref{3:eq:centering}.
\label{4:thm:Tbound}
\end{theorem}
\begin{proof}
See Appendix \ref{Appendix:EstimationT}.
\end{proof}
It follows that if at any Centering Phase $k$ the lower bound in Theorem \ref{4:thm:Tbound} satisfies:
\begin{equation}
    t_k \;\geq \;-\tau -  \left(\left(1+\frac{1}{\alpha\beta}\right)\epsilon_c + \sqrt{2\epsilon_c}\right)||X_{k,L}^{-1}||_F \geq t^*
\end{equation}
the algorithm can be terminated with the optimality gap $\epsilon_g = \frac{N}{t^*}$ guaranteed.

\section{Illustrative Examples}\label{Section:NumericalP}
In this section we illustrate proposed algorithm using examples from the SDPLib dataset \cite{Borchers1999} and randomly generated SDPs.  For testing purposes the algorithm is implemented in Matlab R2016a on a MacBook Pro system with a 3 GHz dual core processor and 16 GB of RAM. The decrease step is carried out on Mosek \cite{mosek} through its Matlab API while the centering step is implemented directly on Matlab using in-house software. The centering epsilon $\epsilon_c$ and global epsilon $\epsilon_g$ are set to $10e-2$  and $10e-3$ for all experiments, respectively.

As a first test bed, we use the \textit{Theta-1} problem from \cite{Borchers1999}, an instance of the Lov\'asz theta number problem, defined  as the solution to the SDP 
\begin{equation}
    \begin{split}
      \vartheta\left( G \right) = \operatorname*{min}_{X \succeq 0}\; &\operatorname{Tr}\left(J^TX\right)\\
        \text{s.t.}\; &\operatorname{Tr}\left(X\right) = 1\\
       &X_{i,j} = 0\quad \left(i,j\right) \in E\\
        \end{split}
    \label{6:Theta}
\end{equation}
where $J$ is the matrix of all ones. The graph for the test problem \textit{Theta-1} has 50 vertices and 103 edges, so the associated SDP \eqref{6:Theta} has $N = 50$ and $M=104$. We solve \eqref{6:Theta} with the proposed algorithm using either the DD cone or the SDD cone in the decrease and centering steps. The results are summarized in Figure  \ref{fig:DDvsSDD}.

\begin{minipage}{0.95\linewidth}
\makebox[\linewidth]{
 \includegraphics[width =0.95 \linewidth]{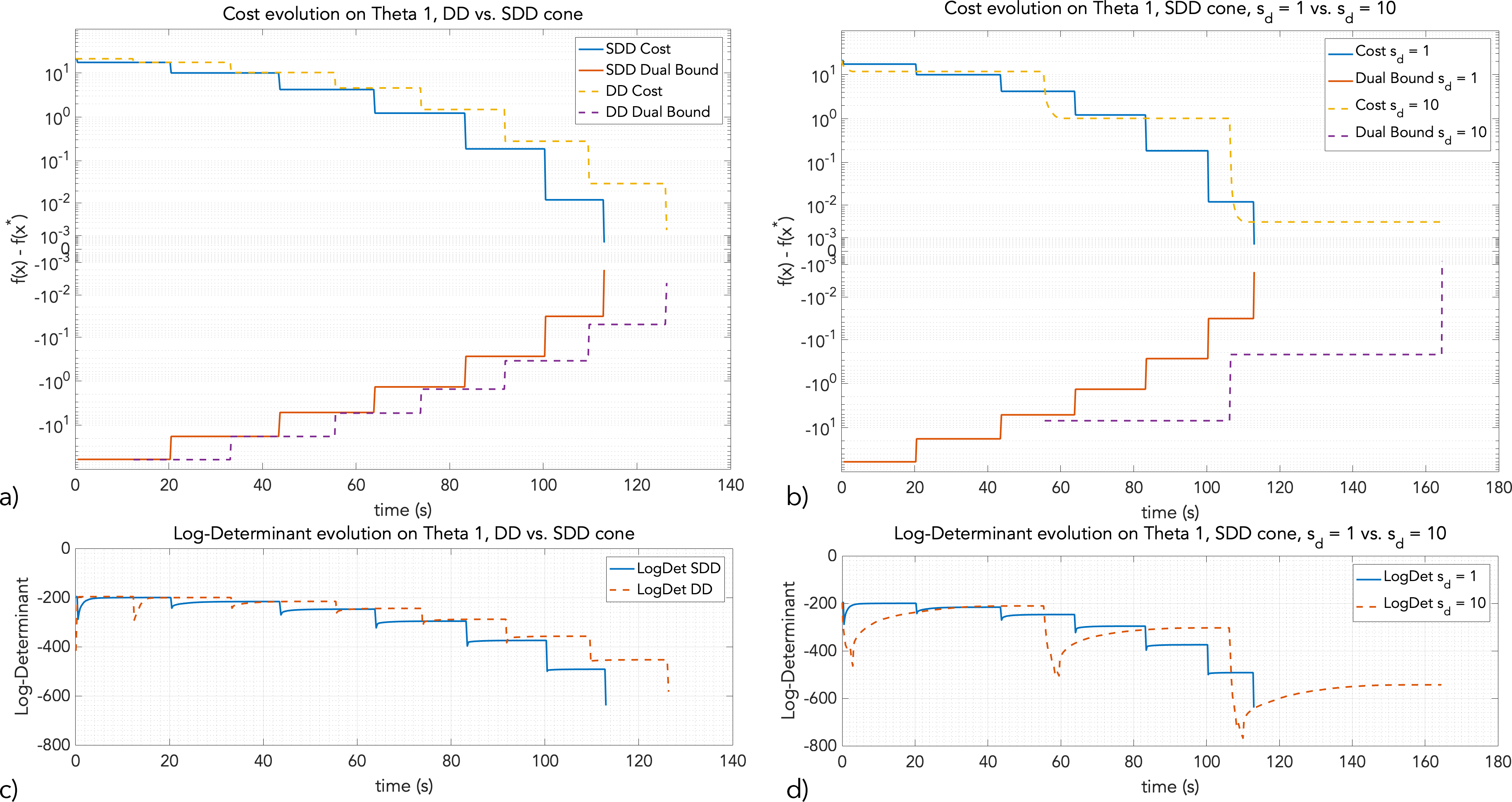}}
\captionof{figure}{$a)$ and $c)$ panels: comparison of algorithm evolution for the DD and SDD cone choices as a function of time for $s_d = 1$. $b)$ and $d)$ panels: comparison of the algorithm evolution using the SDD cone for number of steps $s_d = 1$ and $s_d = 10$.}\label{fig:DDvsSDD}
\end{minipage}

Figure \ref{fig:DDvsSDD} illustrates the evolution of the  proposed algorithm using either the DD cone or the SDD cone in the decrease and centering steps. 
Despite its shorter Centering Phases, the DD approach produces smaller cost improvements at each Decrease Phases and results in more Decrease and Centering Phases and overall a higher runtime. In panels \ref{fig:DDvsSDD}.$b)$ and \ref{fig:DDvsSDD}.$d)$ we show the impact of $s_d$, the number of decrease steps per Phase for the cases  $s_d = 1$ and $s_d = 10$. As expected, choosing a higher $s_d$ results in higher cost decrease per iteration. However, this improvement is offset by the longer Centering Phases, as additional decrease steps tend to result in a worsening of the iterate conditioning, which then requires additional centering steps to bring the iterate back to the central path. In this example, the trade-off between higher cost decrease and longer Centering Phases falls in favor of $s_d = 1$. In Figure \ref{fig:sd} we show a quantitative analysis of this trade-off, where we plot the runtime for both the DD and SDD variant of the algorithm and values of $s_d$ varying from 1 to 10.
\begin{minipage}{\linewidth}
\makebox[\linewidth]{
 \includegraphics[width = 12cm]{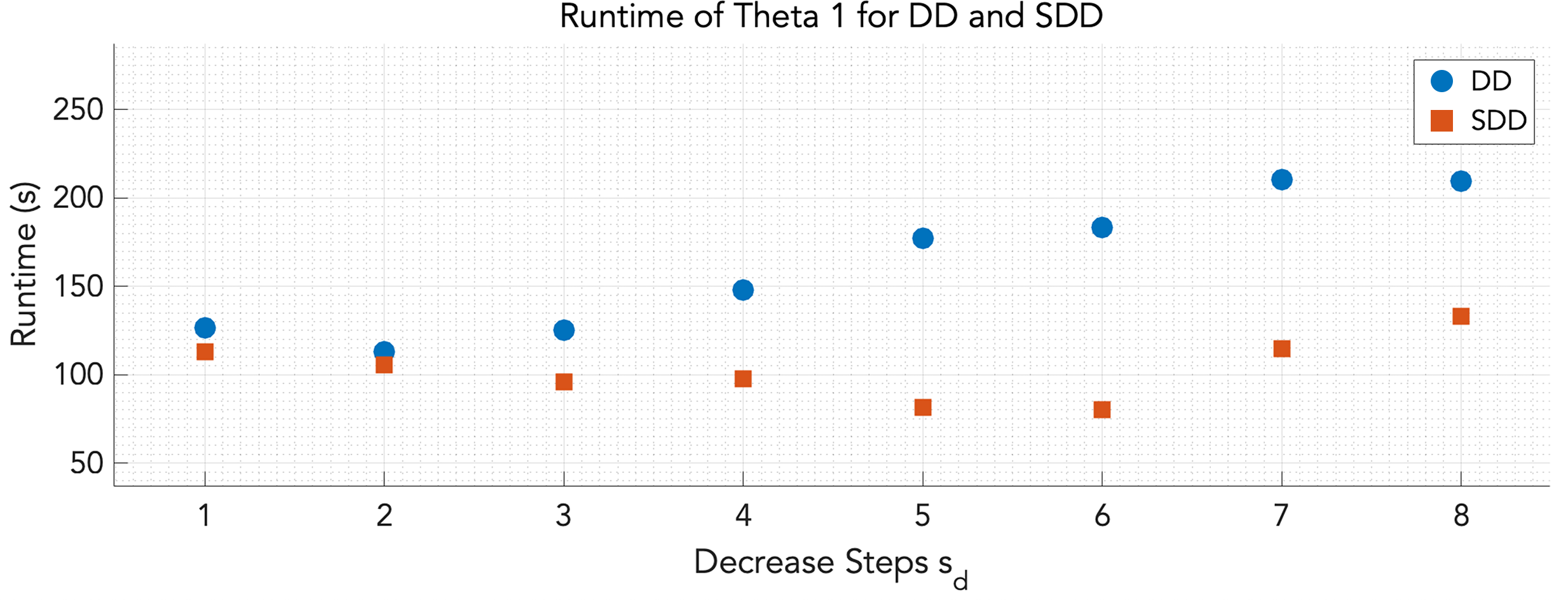}}
\captionof{figure}{Run time of Algorithm \ref{3:Algo:Principal} for DD and SDD cones as a function of the number of decrease steps $s_d$ taken in each Decrease Phase.}\label{fig:sd}
\end{minipage}

Next we test Algorithm 1 on a wider set of SDP problems: from the SDPLib library, we use the theta number instances of the SDPLib \textit{Theta-1} and \textit{Theta-2} and MaxCut problems \textit{mcp100}, \textit{mcp125-1}, \textit{mcp125-2}, \textit{mcp250-1} and \textit{mcp250-2}, which are also SDP relaxations of the well-known NP-hard problem of finding the maximum cut in a graph. We also test the approach on well-behaved randomly generated SDPs. For the proposed algorithm we used the SDD cone for the decrease and centering steps and set $s_d$ to $5$.

\begin{table}[!h]
    \centering
\begin{tabular}{ |c|c|c|c|c|c|c| }
\hline
\multicolumn{3}{|c|}{} & \multicolumn{2}{|c|}{Proposed} & \multicolumn{2}{|c|}{SDD-IBU  \cite{Ahmadi2015}} \\
\hline
Problem & N & M & obj. & time (s) & obj & time (s)\\
\hline \hline
theta1 & 50 & 104 & -23 & 81.69/85.9  & -21 & 37.4  \\
\hline
theta2 & 100 & 498 & -32.88 & 1910.75   & -30.26 & 601.13  \\
\hline
mcp100 & 100 & 100 & -226.16  & 755.5 &  -220.30 & 777.61 \\
\hline
mcp124-1 & 124 & 124 & -141.99 & 1508.3  &   -133.73 & 1338.56  \\
\hline
mcp124-2 & 124 & 124 & -269.88 & 1079.55 &    -264.41 & 1639.61 \\
\hline
mcp250-1 & 250 & 250 & -317.26 & 21259.07 & -288.39 & 5432.37   \\
\hline
mcp250-2 & 250 & 250 & -531.93 & 22074.42 & -485.40 &  6022.99  \\
\hline
RandomSDP-1 & 50 & 50 & 2334.7 & 41.5 & 2625.0 &36.6 \\
\hline
RandomSDP-2 & 50 & 100 & -187.7 & 71.5 & 66.15 & 59.52 \\
\hline
RandomSDP-3 & 100 & 50 & -42.47 & 209.15 & 1.30 & 207.73  \\
\hline
RandomSDP-4 & 100 & 100 & -411.0 & 552.86 &  81.58 & 347.1 \\
\hline
\end{tabular}
    \caption{ Comparison of the proposed algorithm versus SDD with Iterative Basis Update (SDD-IBU) Optimal objective values and execution times for  SDPLib test problems and randomly generated SDPs. In all cases the proposed algorithm achieved the optimal value with an optimality gap $\epsilon_g \leq 0.05$ }
    \label{tab:Table1}
\end{table}

 Table \ref{tab:Table1} summarizes the results of these experiments. As shown there, in all cases the proposed algorithm achieved the optimal value with an optimality gap $\epsilon_g \leq 0.05$, while the SDD with Iterative Basis Update consistently failed to so, in some cases by a large margin.
 
\section{Discussion}\label{Section:Discussion}
In Section \ref{Section:NumericalP} we have shown that Algorithm \ref{3:Algo:Principal}, beyond the theoretical guarantees of Section \ref{Section:Convergence}, also converges in practice and returns $\epsilon_g$-optimal solutions to the original SDP problem \eqref{1:eq:SDP}. However, a salient feature of the results shown in Table \ref{tab:Table1} is that the proposed algorithm is not competitive in terms of runtime with mature IPM solvers, although it should be noted that our implementation of Algorithm \ref{3:Algo:Principal} is a basic proof-of-concept implementation on Matlab. Still, the results shown in Table \ref{tab:Table1} for both the proposed method  and SDD-IBU contrast sharply with the great computational speedups of DD and SDD programs introduced in \cite{Ahmadi2019}. To study this apparent divergence between the efficiency of DD and SDD programs and the perfomance of our method, we analyze the computational complexity of DD and SDD programs and their basis-update extensions and compare them with the computational complexity of IPM solvers for SDPs. 

Recall that the computational complexity of using IPMs to solve  \eqref{1:eq:SDP}  is $\mathcal{O}\left(M N^3 \right.$ \\ $ \left .+ M^2 N^2\right)$ \cite{alizadeh1998}, arising from the Hessian inversion and multiplications needed to compute Newton's step \cite{nesterov1994,Boyd2004}. On the other hand, for  DD and SDD programs, this complexity drops to $\mathcal{O}\left( M^2 N^2\right)$ due to the block-diagonal structure of the  Hessian.  However,  basis-updated DD and SDD programs require computing  the data matrices 
\begin{equation}
    \tilde{C}^T = U_{k-1} C^T U_{k-1}^T,\quad \quad \tilde{A}_i
^T= U_{k-1} A_i^T U_{k-1}^T
\label{6:eq:basischange}
\end{equation}
which entails $M$ matrix products of size $N \times N$ adding   $\mathcal{O}\left(M N^3\right)$ to the original cost. Thus, in principle, the worst case complexity of the algorithm proposed in \cite{Ahmadi2015} is comparable to that of IP methods.  Since the algorithm proposed in this paper inherits the complexity of the basis change algorithm, its worst case asymptotic complexity is also $\mathcal{O}\left(M N^3 + M^2 N^2\right)$.  Note however, that actual complexity for a given problem can be much lower, depending on the number of basis changes required.

\section{Conclusion}\label{Section:Conclusion}
In this work we have developed an algorithm based on the diagonally-dominant (DD) and scaled-diagonally dominant (SDD) SDP relaxations developed in \cite{Ahmadi2019,Ahmadi2015} to solve SDPs to  $\epsilon_g$-global optimality in a polynomially bounded number of iterations. We note that, as presented, the proposed algorithm is not competitive  vis-a-vis  mature, commercially available IPM solvers, although it should be noted that our implementation of Algorithm \ref{3:Algo:Principal} is a basic proof-of-concept implementation on Matlab. Nevertheless, we believe that the proposed approach is valuable for two reasons.  First, from a theoretical standpoint, our work shows that it is indeed possible to solve an SPD to $\epsilon$-optimality by solving  DD and SDD conic problems, giving a positive answer to a question left open in  \cite{Ahmadi2015,Ahmadi2019}.  From a 
 practical side, Algorithm \ref{3:Algo:Principal} can be used as a blueprint for the development of efficient algorithms based on the notion of combining a sequence of simple conic problems with a \textit{decrease-and-center} iterative methodology. For instance, other conic sets could be used in an analogous manner in Algorithm \ref{3:Algo:Principal}. A straightforward choice would be any cone in the factor-width cone family $\mathcal{FW}_k$ \cite{Boman2005,gouveia2019sums}, of which the DD and SDD cones are members of for $k=1,2$, respectively. Beyond those, any full-dimensional cone satisfying Lemma \ref{4:lemma:sphere} and whose barrier is self-concordant and can bound an affine function on the log-determinant as in Corollary \ref{4:coro:Id} could be a suitable choice for Algorithm \ref{3:Algo:Principal}.
Other interesting alternatives are structured cones with favorable numerical properties. While usually not full-dimensional, the structure of these cones can be leveraged nonetheless to increase the performance of the optimization. In that regard, we have investigated the use of intersections of the DD and SDD cones with the cone of symmetric band matrices with band size $k$.  Elements within this conic intersection can be parameterized by $\mathcal{O}\left(Nk\right)$ variables, and the barrier function for this set can be defined by the original barrier functions of the DD/SDD sets defined in \eqref{4:eq:logbarrs}, limiting the summation only to the terms satisfying $j-i\leq k$. Using these simplified barriers greatly reduces the cost of the DD and SDD problem and brings the cubic complexity of the Cholesky factorization and the basis change with respect to $N$ down to linear complexity.  Finally, we are investigating the possibility of combining the proposed algorithm with the cone decompositions proposed in \cite{MILLER2022}.

\section*{Acknowledgments}
This work was partially supported by NSF grants  IIS--1814631, ECCS—1808381 and CNS--2038493, and AFOSR grant FA9550-19-1-0005.
\bibliographystyle{siamplain}
\bibliography{references}

\appendix 

\section{Proofs of Section \ref{Subsec:Centering}}\label{Appendix:ProofsCentering}

\subsection{Lemma \ref{4:lemma:selfconcor}}\label{App:PC:selfconcor} To prove the Lemma we need the following additional result:

\begin{lemma} Consider the set  $\mathcal{M} = \{M_{1,2}, M_{1,3},\dots,M_{N-1,N}\}$ consisting of $N \choose 2$ $2 \times 2$ matrices. Let
 $Y = \sum  \Psi_{i,j}(M_{i,j}) \doteq \Psi(\mathcal{M}) \in SDD_N$. Assume that $Y \in SDD_N$,  $N$ is even and $  \phi_{SDD}(\mathcal{M})$ is finite. Then there exists a set of $N-1$ matrices $Z_k\succ 0$ such that $Y = \sum_k^{N-1} Z_k$ and $\phi_{SDD}(\mathcal{M}) = \sum \log(|Z_k|)$.
\label{app:lemma:tournament}
\end{lemma}
\begin{proof} Our goal is to decompose the 
 set $\mathcal{M}$ into $N-1$ disjoint subsets $\mathcal{M}_k$ such that $Z_k = \Psi(\mathcal{M}_k)$. To do so, define a graph $G = (V,E)$ on $N$ vertices and  a set $\mathcal{S}=\{S_1,\dots,S_{N-1}\}$ where each element $S_l$ of the set $\mathcal{S}$ is a subset of the edges of $G$, i.e. $S_l \subseteq E$ for all $l=1,\dots,N-1$. Assume furthermore that $G$ is a complete graph $K_N$, i.e. all the possible edges of $G$ are contained in $E$.
An edge-coloring of a graph $G$ is an assignment of labels (colors) to the edges of the $G$ so that no two incident edges have the same label. Given an edge-coloring of $G$, define $S_k$ as the set of all edges corresponding to the label $k$. For the case of complete graphs, it is well-known that $K_N$ for even $N$ is edge-colorable with $N-1$ colors \cite{Behzad1967}, i.e. there exist $N-1$ sets of edges $S_k$ such that each set induces a perfect matching on $K_N$ and each edge of $K_N$ appears in exactly one set of $\mathcal{S}$. In other words, $S_1 \cup \dots,\cup S_{N-1} = E$  and $S_k \cap S_l = \emptyset$ for any $k\neq l$.
Assume that $\mathcal{S}$ is the edge set produced by an edge-coloring of $K_N$. Assume also that $M_{i,j} \in \mathcal{M}_k$ if and only if $\{i,j\} \in S_k$, i.e. the sets $\mathcal{M}_k$ contain the matrices $M_{i,j}$ labeled by the edges $\{i,j\}$ present in each coloring $S_k$. Then defining $Z_k = \Psi(\mathcal{M}_k)$ implies that $Y = \sum_k^{N-1} Z_k$, as each block $M_{i,j}$ contributes to exactly one of the matrices $Z_k$. Since
each set of indices $S_k$ is a perfect matching, there must exist a permutation matrix $P_k$ such that $P_k Z_k P_k^T$ is a block-diagonal matrix whose blocks are the $2\times 2$ $M_{i,j}$ matrices in $\mathcal{M}_k$. Furthermore, finitiness of $\phi_{SDD}(\mathcal{M})$ and the membership of $Y$ in $SDD_N$ guarantee that $Z_k\succ 0$, leading to:
\begin{equation}
\log(|Z_k|) = \log(|P_k Z_k P_k^T|) = \sum_{\{i,j\} \in S_k}\;\log(|M_{i,j}|)
\end{equation}
Summing over all the indices $k$ leads to the final result:
\begin{equation}
    \sum_k^{N-1} \log(|Z_k|) =  \sum_k^{N-1} \sum_{\{i,j\} \in S_k}\;\log(|M_{i,j}|) = \sum_{i,j>i}^N\;\log(|M_{i,j}|) = \phi_{SDD}(\mathcal{M})
\end{equation}
\end{proof}
\noindent Using this result, we can now prove Lemma \ref{4:lemma:selfconcor}
\begin{proof}
\noindent \textbf{Property (a):}
For the DD case,   each term of the form $-\log\left(M_{i,j}(1,1)^2 \right.$-$\left . M_{i,j}(1,2)^2\right)$ can be written as $-\log\left(M_{i,j}(1,1)+ M_{i,j}(1,2)\right)$  \\ $-\log\left(M_{i,j}(1,1) \right.$ -$\left. M_{i,j}(1,2)\right)$. As the negative logarithm of an affine function is self-concordant with respect to the arguments of said function, and self-concordance is preserved under summation \cite{Boyd2004}, $-\phi_{DD}(\mathcal{M})$ is self-concordant. For the barrier of the SDD set, each term in the summation can be expressed as $-\log\left(|M_{i,j}|\right)$. Since  the negative log-determinant is self-concordant with respect to its argument,  $-\phi_{SDD}(\mathcal{M})$ is also self-concordant.

\noindent \textbf{Property (b):} We prove the result by showing that the inequality above also holds individually for each term in the summations in \eqref{4:eq:logbarrs}. For clarity, write $M_{ij}$ as $\begin{pmatrix}a&c\\c&b\end{pmatrix}$. Then this term-by-term inequality implies that:
\begin{equation*}
    -\frac{1}{2}\log(a^2-c^2) -\frac{1}{2}\log(b^2-c^2) \geq -\log(ab -c^2)
\end{equation*}
Multiplying by $-2$ and applying an exponential to both sides leads to:
\begin{equation*}
    (a^2-c^2)(b^2-c^2) \leq (ab-c^2)^2 \implies a^2c + b^2c - 2abc^2 \geq 0 \implies (a-b)^2c^2\geq 0
\end{equation*}
which holds trivially, finishing the proof.

\noindent  \textbf{Property (c):} In the following we only consider the case for which $N$ is even, as any SDP like \eqref{1:eq:SDP} whose matrix variables are of odd dimension $N$ can be embedded in an SDP of even dimension $N+1$ by adding trivial constraints on the $N+1$'th row and column of the embedding matrix variable, for which case the following results hold. As such, assuming $N$ is even and using the results from Lemma \ref{app:lemma:tournament}, have that:
\begin{equation}
    \begin{split}
        Y &= \sum_k^{N-1} Z_k\\
        \frac{1}{N-1}|Y|^{1/N} &= \frac{1}{N-1}\mid\sum_k^{N-1} Z_k\mid^{1/N}\\
         \frac{1}{N-1}|Y|^{1/N} &\geq \frac{1}{N-1}\sum_k^{N-1}| Z_k|^{1/N}\\
    \end{split}
\end{equation}
where the inequality comes from Minkowski's determinant inequality. Applying the logarithm to both sides and using its concavity yields:
\begin{equation}
        \begin{split}
      &   \log\left(\frac{1}{N-1}|Y|^{1/N} \right)\geq \log\left(\frac{1}{N-1}\sum_k^{N-1}| Z_k|^{1/N}\right)\\
     &    \frac{1}{N} \log(|Y|) - \log\left(N-1\right)\geq \log\left(\frac{1}{N-1}\sum_k^{N-1}| Z_k|^{1/N}\right) \\
      &   \frac{1}{N} \log(|Y|) - \log\left(N-1\right)\geq \frac{1}{N-1}\sum_k^{N-1}\frac{1}{N}\log(| Z_k|) \\
     &   h(Y) =  \left(N-1\right) \log(|Y|) - N\left(N-1\right) \log\left(N-1\right)\geq \sum_k^{N-1}\log(| Z_k|) \\
     &    h(Y)\geq \phi_{SDD}(\mathcal{M}) \; \implies \;-\phi_{SDD}(\mathcal{M})\geq-h(Y)
    \end{split}
\end{equation}
where the last inequality follows from Lemma \ref{app:lemma:tournament}.\\
\\
\end{proof}

\subsection{Proof of Lemma \ref{4:lemma:SCdamped}}\label{App:PC:SCdamped}

\begin{proof}
Take $\mathcal{M}_+$ to be the result of taking a Newton step from $\mathcal{M}_0$, where $\mathcal{M}_0$ is the set where all its elements are of the form $M_{i,j} = \frac{1}{N-1} I_{2\times 2}$, i.e. $Y_0 = \Psi(\mathcal{M}_0) = I$, and by construction is always a feasible solution of \eqref{3:eq:centering}. Then by \eqref{4:eq:sc1} we have that:
\begin{equation*}
    -\phi(\mathcal{M}_+) \leq -\phi(\mathcal{M}_0)  - \alpha \beta \frac{ \lambda^2(x)}{1+ \lambda(x)}  \leq -\phi(\mathcal{M}_0)- \xi
\end{equation*}
By Corollary \ref{4:coro:Id}, the right hand side is equivalent to $N(N-1)\log(N-1) - \xi$. Taking the $\mathcal{M}$ that optimizes problem \eqref{3:eq:centering} as $\mathcal{M}*$, the left hand side can be lower bounded using also Corollary \ref{4:coro:Id} as:
\begin{equation*}
\begin{split}
        &   -\phi(\mathcal{M}_+) \geq  -h(\Psi(\mathcal{M}_*)) = -\left(N-1\right) \log\left(|\Psi(\mathcal{M}_*)|\right) + N \left(N-1\right)\log\left(N-1\right)\\
           &=-\left(N-1\right)\left(\log\left(|X_l|\right)-\log\left(|X_{l-1}|\right)\right)  + N \left(N-1\right)\log\left(N-1\right)
\end{split}
\end{equation*}
where we have used the fact that in problem \eqref{3:eq:centering} we have that $U_{l-1}^T Y_* U_{l-1} = X_l$ and the properties of the log-determinant. Combining the last two equations leads to:
\begin{equation*}
\begin{split}
       -\left(N-1\right)\left(\log\left(|X_l|\right)-\log\left(|X_{l-1}|\right)\right)  + N \left(N-1\right)\log\left(N-1\right) &\leq N(N-1)\log(N-1) \\ - \xi \Rightarrow 
            -(N-1)\log\left(|X_l|\right) \leq -(N-1)\log\left(|X_{l-1}|\right) - \xi
\end{split}
\end{equation*}
which completes the proof.
\end{proof}
\subsection{Lemma \ref{4:lemma:SCquad}}\label{App:PC:SCquad} The following
additional Lemma is needed for the proof:
\begin{lemma}
The Newton decrements $\lambda_{\phi}(\mathcal{M})$ and $\lambda_{h}(\Psi(\mathcal{M}))$ of problems \eqref{3:eq:centering} and \eqref{4:eq:PSDcent}, respectively, evaluated at $\mathcal{M}_0$ satisfy:
\begin{equation}
   \lambda_{\phi}(\mathcal{M}_0) \geq \frac{\lambda_{h}(\Psi(\mathcal{M}_0))}{\sqrt{N-1}}
\end{equation}
\label{4:lemma:lambdas}
\end{lemma}
\begin{proof}
For simplicity, in the following we use a vectorized notation for the sets $\mathcal{M}$, where we collapse the whole set $\mathcal{M}$ onto a vector $m\in\mathbb{R}^{3{N \choose 2}}$ by stacking the vectorization of all $M_{i,j}$ matrices as $m = [m_{1,2},m_{1,3},\dots,m_{N-1,N}]^T$, where $m_{i,j} = [M_{i,j}(1,1),M_{i,j}(2,2),M_{i,j}(1,2)]$. Using this vectorized notation, the Newton decrement for problem \eqref{3:eq:centering} evaluated at $m$ is bounded below by:
\begin{equation}
    \lambda_{\phi}(m) \geq  \frac{v^T \nabla_{\phi}(m)}{\left( -v^T \nabla^2_{\phi}(m) v \right)^{1/2}}
\end{equation}
for any feasible $v$, with equality for the Newton step $v = m_{nt}$ \cite{Boyd2004}. Take the Newton step $\Delta_{nt}$ for problem \eqref{4:eq:PSDcent} evaluated at $Y=I$ and define $\tilde{\mathcal{M}}_{nt} = \{\dots,\tilde{M}_{i,j}\dots\}$ with $\tilde{M}_{i,j}(1,1) = \Delta_{nt}(i,i)/(N-1)$, $\tilde{M}_{i,j}(2,2) = \Delta_{nt}(j,j)/(N-1)$ and $\tilde{M}_{i,j}(1,2) = \tilde{M}_{i,j}(2,1) = \Delta_{nt}(i,j)$, from which follows that $\Psi(\tilde{\mathcal{M}}_{nt}) = \Delta_{nt}$. Evaluating the above expression at $m = m_0$ taking $v = \tilde{m}_{nt}$, the vectorization of $\tilde{\mathcal{M}}_{nt}$, leads to:
\begin{equation}
\begin{split}
      &  \lambda_{\phi}(m_0) \geq  \frac{\tilde{m}_{nt}^T \nabla_{\phi}(m)}{\left( -\tilde{m}_{nt}^T \nabla^2_{\phi}(m_0) \tilde{m}_{nt} \right)^{1/2}} \\ & = \frac{\operatorname{Tr}(\Delta_{nt})}{\left( \sum_i^N (N-1)\left( \Delta_{nt}(i,i)\right)^2 + \sum_{i,j<i} 2(N-1)^2\left( \Delta_{nt}(i,j)\right)^2 \right)^{1/2}} \\
         &\geq \frac{\operatorname{Tr}(\Delta_{nt})}{\left( (N-1)^2 \operatorname{Tr}(\Delta_{nt}^T\Delta_{nt}) \right)^{1/2}}= \frac{\operatorname{Tr}(\nabla_h(Y)^T\Delta_{nt})}{\sqrt{N-1}\left( -\Delta_{nt}^T \nabla^2_h(Y) \Delta_{nt} \right)^{1/2}} = \frac{\lambda_{h}(\Psi(\mathcal{M}_0))}{\sqrt{N-1}}
\end{split}
\end{equation}
where we have used the gradient and Hessian evaluations presented in Section \ref{Appendix:GradsHessians}.
\end{proof}

Using the results above, we can proceed to the proof of Lemma \ref{4:lemma:SCquad}:

\begin{proof}
From  the hypothesis and Lemma \ref{4:lemma:lambdas}, we have that $\lambda_h(\Psi(\mathcal{M}_0))\leq \eta$. Since $\eta < 0.68$ by construction,  \eqref{4:eq:sc3} yields the bound on the optimal value of problem \eqref{4:eq:PSDcent}:
\begin{equation*}
    -h(\Psi(\mathcal{M}_0)) \geq -h(Y_*) \geq  -h(\Psi(\mathcal{M}_0)) - \lambda_h^2(\Psi(\mathcal{M}_0)) \geq -h(\Psi(\mathcal{M}_0)) -(N-1)\lambda^2_{\phi}(\mathcal{M}_0)
\end{equation*}
Applying equation \eqref{4:eq:sc2} and Corollary \eqref{4:coro:Id} to problem \eqref{3:eq:centering} leads to:
\begin{equation}
    -\phi(\mathcal{M}_0) - \alpha \lambda_\phi^2(\mathcal{M}_0)\geq -\phi(\mathcal{M}_+) \geq -h(Y_l) 
    \label{4:eq:scproofeq1}
\end{equation}
Combining the last two equations leads to:
\begin{equation*}
    \begin{split}
   &   h(\Psi(\mathcal{M}_0)) + (N-1)\lambda^2_{\phi}(\mathcal{M}_0) -\phi(\mathcal{M}_0) - \alpha \lambda_\phi^2(\mathcal{M}_0) \geq -h(Y_l) +   h(Y_*)\\
  &   \Rightarrow (N-1-\alpha) \lambda_\phi^2(\mathcal{M}_0) \geq (N-1)\,\left(\log\left(|X_*|\right) -\log\left(|X_l|\right)\right)   
      \end{split}
\end{equation*}
Setting $g_l =  (N-1)\,\left(\log\left(|X_*|\right) -\log\left(|X_l|\right)\right) $ yields $g_l \leq  (N-1-\alpha) \lambda_\phi^2(\mathcal{M}_0) \leq \frac{N-1-\alpha}{N-1} \eta^2$.\\
\\
If $ \lambda_\phi^2 \leq \frac{g_{l-1}}{N-1}$, then $g_l$ is also upper-bounded by $\frac{N-1-\alpha}{N-1} g_{l-1} \geq g_l$. Otherwise, if $ \lambda_\phi^2 >  \frac{g_{l-1}}{N-1}$, combining equation \eqref{4:eq:scproofeq1} with the Lemma's assumptions leads to:
\begin{equation*}
    \begin{split}
      &  -\phi(\mathcal{M}_0) - \alpha \lambda_\phi^2(\mathcal{M}_0) + g_{l-1}  \geq -h(Y_l) + (N-1)\,\left(\log\left(|X_*|\right) -\log\left(|X_{l-1}|\right)\right)\\
      &     g_{l-1} - \alpha \lambda_\phi^2(\mathcal{M}_0)  \geq -h(Y_l) +\phi(\mathcal{M}_0) + (N-1)\,\left(\log\left(|X_*|\right) -\log\left(|X_{l-1}|\right)\right)\\
       &  g_{l-1} - \alpha \lambda_\phi^2(\mathcal{M}_0)  \geq (N-1)\,\left(\log\left(|X_{l-1}|\right) -\log\left(|X_{l}|\right)\right)+\\& \phantom{g_{l-1} - \alpha \lambda_\phi^2(\mathcal{M}_0)} + (N-1)\,\left(\log\left(|X_*|\right) -\log\left(|X_{l-1}|\right)\right)\\
       &   \frac{N-1-\alpha}{N-1} g_{l-1} \geq g_{l-1} - \alpha \lambda_\phi^2(\mathcal{M}_0)  \geq  (N-1)\,\left(\log\left(|X_*|\right) -\log\left(|X_{l}|\right)\right) = g_l
    \end{split}
\end{equation*}
From which follows that $\frac{N-1-\alpha}{N-1} g_{l-1}\geq g_l$ for any value of $\lambda_\phi^2$, finishing the proof.
\end{proof}

\section{Proofs of Section \ref{Subsec:Conv:Dec}}\label{Appendix:ProofsDecrease}

\subsection{Lemma \ref{4:lemma:decrease}}\label{App:PD:decrease} First we introduce two results from \cite{Boman2005} that will be used in the proof.
\begin{definition*}[adapted from \cite{Boman2005}]
A symmetric matrix $A$ is said to be an H-matrix if the matrix $M(A)$  defined by:
\begin{equation}
    (M(A))_{ij} = 
    \begin{cases}
   |a_{ij}|,      & \quad i=j\\
    -|a_{ij}|,      & \quad i \neq j\\
  \end{cases}
\end{equation}
is positive semidefinite.
\end{definition*}

\begin{theorem}\label{teo:boman}[Adapted from (8) from \cite{Boman2005}]
A symmetric matrix $A$ is an H-matrix if and only if $A$ is scaled diagonally dominant.
\end{theorem}

To prove  Lemma  \ref{4:lemma:decrease} we need the following result:
\begin{lemma} \label{4:lemma:sphere}
If $||Y-I||_F^2 \leq 1$, then $Y$ is scaled diagonally-dominant. Furthermore, if $||Y-I||_F^2 \leq \frac{2}{N+1}$, $Y$ is diagonally-dominant.
\end{lemma}
\begin{proof}
We start with the proof for the SDD case. First note the following fact:
\begin{equation}
    ||Y-I||_F^2 \leq 1 \implies Y\succeq 0
\end{equation}
which follows from the fact that the Frobenius norm can also be expressed as a function of the eigenvalues of $Y$ as in $||Y-I||_F^2 = ||\Sigma_Y -I||_F^2 \leq 1$, which implies that no eigenvalue of $Y$ can be lower than $0$, thus proving positive-semidefiniteness. A consequence of this is that the diagonal values of $Y$ must satisfy $Y_{ii}\geq0$.

Note also that the Frobenius norm bound on $Y-I$ also extends to its comparison matrix, i.e. $||M(Y-I)||_F^2 \leq 1$, as the Frobenius norm acts element-wise on its argument and is not affected by the signs of the entries. Due to the non-negativity of the diagonal entries $Y_{ii}$, we have that  $||M(Y)-I||_F^2= ||M(Y-I)||_F^2\leq 1$ and that $M(Y)$ is positive-semidefinite. Hence, $Y$ is an $H$-matrix and,  from Theorem \ref{teo:boman}, SDD.\\
\\
For the DD case, assume by contradiction that there exists a $Y \notin DD$  satisfying $||Y-I||_F^2 \leq \frac{2}{N+1}$. Since  $Y \notin DD$ there exists an index $i$ such that:
\begin{equation}
    |Y_{i,i}| < \sum_{j\neq i} |Y_{i,j}|
\end{equation}
Lower-bounding the Frobenius norm by evaluating only the i'th row and column of $Y$ yields:
\begin{equation}
\begin{split}
      ||Y-I||_F^2 &\geq \left(Y_{i,i}-1\right)^2 + 2\sum_{j\neq i} \left(Y_{i,j}\right)^2\\
      &\geq \left(Y_{i,i}-1\right)^2 +\frac{2}{N-1}\left(\sum_{j\neq i} |Y_{i,j}|\right)^2 \\
      &> \left(Y_{i,i}-1\right)^2 +\frac{2}{N-1}Y_{i,i}^2 = \frac{N+1}{N-1}Y_{i,i}^2 +1 - 2 Y_{i,i}\\
      &\geq \frac{2}{N+1}
\end{split}
\end{equation}
where the first inequality comes from the partial evaluation of the Frobenius norm, the second from the inequality $|z|_2^2\geq \frac{|z|^2_1}{N}$, the third from the non-diagonal dominance of $Y$ and the fourth from the fact that $\frac{2}{N+1}$ is the minimum value of the convex quadratic polynomial $\frac{N+1}{N-1}z^2 - 2 z +1$. The chain of inequalities leads to $||Y-I||_F^2>\frac{2}{N+1}$ which contradicts the initial assumption on the Frobenius norm.  
\end{proof}

\noindent \emph{Proof of Lemma  \ref{4:lemma:decrease}.}
The dual of problem  \eqref{4:eq:BarrierPr} is:
\begin{equation}
    \begin{split}
        Z_t, y_t = \operatorname*{argmin}_{Z\succeq 0,y} &\quad -b^T y\; -\frac{1}{t}\log\left(|Z|\right)\\
        \text{s.t.}& \quad Z = C  - \sum_{i=1}^M y_i \,A_i
    \end{split}
    \label{4:eq:BarrierDual}
\end{equation}
From duality, it follows that $X_t = \frac{1}{t} Z_t^{-1}$ and $y_t = \gamma_{i,t}$ \cite{Boyd2004}, where $\gamma_{i,t}$ are the dual variables of the linear constraints of \eqref{4:eq:BarrierPr}. Furthermore, 
from Proposition 6 from \cite{Drummond2002} it follows that the dual variables $\gamma_{i,t}$ evaluated at $X = X_t$ are bounded and thus there must exist a finite quantity $\Theta_t$ such that $\sum_i \gamma_{i,t}^2 \leq \Theta_t$,. Finally, from the finitiness of $\Theta_t$, there must exist a finite $\Theta$ such that $\Theta \geq \Theta_t$ on the finite interval $t_0\geq t \geq t^*$. The KKT conditions of \eqref{4:eq:BarrierPr} yield:
 \[  X_{k-1}^{-1} =  X_{t_{k-1}}^{-1} = t_{k-1} \left(C - \sum \gamma_i A_i\right)\]
Hence, orthogonality of the data matrices $C, A_i$ implies:
\begin{equation}
    ||X_{t_{k-1}}^{-1}||_F^2  = t_{k-1}^2  \left(1 + \sum \gamma_i^2\right) \leq  t_{k-1}^2  \left(1 + \Theta \right)
    \label{eq:dualKKT}
\end{equation}
 By construction $\hat{Y}_k = I - \sqrt{\Phi}\frac{Q}{||Q||_F}$  satisfies $||\hat{Y}_k-I||_F^2 = \Phi$, ensuring by Lemma \ref{4:lemma:sphere} that $\hat{Y}$ is contained in the appropriate DD or SDD cone. Moreover,  the linear constraints evaluated at $\hat{X}_k =U_{k-1}^T \hat{Y}_k U_{k-1}$ satisfy:
\[\begin{aligned}
   &     \operatorname{Tr}\left(A_i^T \hat{X}_k\right)  = \operatorname{Tr}\left(A_i^T U_{k-1}^T \hat{Y} U_{k-1}\right) \\ &= \operatorname{Tr}\left(A_i^T X_{k-1}\right)- \frac{\sqrt{\Phi}}{||Q||_F}\operatorname{Tr}\left(A_i^T U_{k-1}^T Q U_{k-1}\right)\\
        &= \operatorname{Tr}\left(A_i^T X_{k-1}\right)- \frac{\sqrt{\Phi}}{||Q||_F}\operatorname{Tr}\left(A_i^T C\right)
        = \operatorname{Tr}\left(A_i^T X_{k-1}\right) = b_i
\end{aligned}
\]
due to the orthogonality between $A_i$ and $C$  (Assumption \ref{4:ass:1}). Thus $\hat{Y}_k$ is a feasible solution of \eqref{2:eq:IBUDDSDD} with associated cost: \begin{equation} \label{4:eq:costQ1} \begin{aligned}
   & \hat{\mathcal{C}}_k  = \operatorname{Tr}\left(C^TU_{k-1}^T \hat{Y} U_{k-1}\right)  =  \operatorname{Tr}\left(C^T X_{k-1}\right) -\frac{\sqrt{\Phi}}{||Q||_F}
    \end{aligned}
\end{equation}
since  $||C||_F^2 = 1$ from Assumption \ref{4:ass:1}. Using the properties of the Frobenius norm, the norm of $Q$ can be bounded as $||Q||_F = ||U_{k-1}^{-T} C U_{k-1}^{-1}||_F \leq ||C||_F||X_{k-1}^{-1}||_F$. The matrix $X_{k-1}$ is the optimizer of the Centering Phase  at stage $k-1$ and thus by Lemma \ref{4:lemma:CentEquiv} satisfies $X_{k-1} = X_{t_{k-1}}$.  Combining this with  \eqref{eq:dualKKT}
yields:
\begin{equation}
\begin{aligned}
   &  \operatorname{Tr}\left(C^T X_k\right)  \leq \operatorname{Tr}\left(C^T \hat{X}_k\right) = \operatorname{Tr}\left(C^T X_{k-1}\right) -\frac{\sqrt{\Phi}}{||Q||_F} \\
 &    \leq\operatorname{Tr}\left(C^T X_{k-1}\right)  - \frac{\sqrt{\Phi}}{||C||_F||X_{k-1}^{-1}||_F}  \\ &
     \leq \operatorname{Tr}\left(C^T X_{k-1}\right) - \frac{\sqrt{\Phi}}{t_{k-1}\sqrt{1+\Theta}} <\operatorname{Tr}\left(C^T X_{k-1}\right)   
\end{aligned}
    \label{4:eq:costlb}
\end{equation}
\qed

\subsection{Strict decreasing property of $C^TX_{t}$  in \eqref{4:eq:BarrierPr}}\label{App:PD:alb}
We will prove the following more general result:
\begin{lemma} Consider the following optimization
\begin{equation}
    \operatorname*{minimize}_{x\in\mathcal{G}}\;a(x) + \lambda b(x)
    \label{4:eq:alb}
\end{equation}
where $a(x)$ and $b(x)$ are convex and strongly convex functions of $x$, respectively, $\mathcal{G}$ is a convex set and $\lambda \geq 0$ is a tradeoff parameter. Take $x_1$ and $x_2$ to be the optimizers of \eqref{4:eq:alb} corresponding to  $\lambda=\lambda_1$ and  $\lambda = \lambda_2$. 
If $\lambda_1 > \lambda_2$, then $a(x_1) > a(x_2)$.
\label{4:lemma:alb}
\end{lemma}

\begin{proof}
By the strong convexity of the objective function, we have that:
\begin{equation}
    \begin{split}
        a(x_1) + \lambda_1 b(x_1) < a(x_2) + \lambda_1 b(x_2)\\
        a(x_2) + \lambda_2 b(x_2) < a(x_1) + \lambda_2 b(x_1)
    \end{split}
\end{equation}
Combining these two equations leads to:
\begin{equation}
    (\lambda_1-\lambda_2) b(x_1) < (\lambda_1-\lambda_2) b(x_2)
\end{equation}
which implies that $b(x_2)>b(x_1)$. Finally, this last result leads to:
\begin{equation}
   a(x_1)- a(x_2) > \lambda_2 (b(x_2)-b(x_1)) \geq 0
\end{equation}
which finishes the proof.
\end{proof}

\subsection{Lemma \ref{4:lemma:gt}}\label{App:PD:gt}

\begin{proof}
Multiplying the KKT stationarity of problem \eqref{4:eq:BarrierPr} in Lemma \ref{4:lemma:CentEquiv} by $X_t$ on the right and taking the trace leads to:
\begin{equation}
\begin{split}
 &  C -\frac{1}{t}X^{-1}_t-  \sum_i^M \gamma_{it} A_i = 0 \\ & \implies \operatorname{Tr}\left( \left(C -\frac{1}{t}X^{-1}_t- \sum_i^M \gamma_{it} A_i\right)^T X_t \right) = 0 \\
&   \implies \operatorname{Tr}\left( C^T X_t -\frac{1}{t}I - \sum_i^M \gamma_{it} A_i^T X_t \right) = 0\\
&   \implies  \operatorname{Tr}\left( C^T X_t\right) - \frac{N}{t} - \gamma^T b = 0
   \implies \operatorname{Tr}\left( C^T X_t\right)  = b^T y_t + \frac{N}{t}
\end{split}
\end{equation}
where we have used the fact that the dual variables $\gamma_{it}$ are equivalent to $y_t$ in \eqref{4:eq:BarrierDual}. Applying Lemma \ref{4:lemma:alb} to
  \eqref{4:eq:BarrierDual} shows that
 $b^T y_t$ is monotonically increasing. Thus, for any $t\geq t_0$
\begin{equation}
    \operatorname{Tr}\left( C^T X_t\right) = b^T y_t + \frac{N}{t} \geq b^T y_{t_0} + \frac{N}{t} = \operatorname{Tr}\left( C^T X_{t_0}\right) - \frac{N}{t_0} + \frac{N}{t} = g_0(t)
\end{equation}
which finishes the proof.
\end{proof}

\section{Proof of Theorem \ref{4:thm:Tbound}}\label{Appendix:EstimationT}
\begin{proof}
We start the proof by stating the stationarity KKT condition of \eqref{4:eq:BarrierPr} at $X_t$ and the stationarity KKT condition of problem \eqref{4:eq:PSDcent2} at $X_L$:
\begin{equation}
    \begin{split}
        -X^{-1}_t &=  -t C + t\sum \gamma_i A_i \\
      E   -X_L^{-1} &= \tau C + \sum \beta_i A_i 
    \end{split}
    \label{App:eq:tboundKKT}
\end{equation}
where $E$ is a residual of the KKT stationarity condition for $X_L$, due to its suboptimality. From these equation it follows that:
\begin{equation}
\begin{split}
   &  ||X_L^{-1} -X_t^{-1} - E||_F^2 = ||\left(-\tau-t\right) C + \sum \left(t \gamma_i -\beta_i \right) A_i||_F^2\\
    & = \left(-\tau-t\right)^2 ||C||_F^2+ \sum \left(t \gamma_i -\beta_i\right)^2 ||A_i||_F^2
      \geq \left(-\tau-t\right)^2 
\end{split}
\end{equation}
where the equality comes from the orthogonality between $C$ and $A_i$ and the inequality from the non-negativity of the terms dependent on $||A_i||_F^2$ and the unit norm of $C$. Taking the square root on both sides of the inequality and applying the triangle inequality on the left hand side leads to:
\begin{equation}
     || E||_F +  ||X_L^{-1} -X_t^{-1}||_F \geq |-\tau-t|
\end{equation}
The second Frobenius norm can be bounded as follows:
\begin{equation}
\begin{split}
    ||X_L^{-1} -X_t^{-1}||_F 
    &= ||X_L^{-1/2}\left(I -X_L^{1/2}X_t^{-1}X_L^{1/2}\right)X_L^{-1/2}||_F\\
     &\leq ||I -X_L^{1/2}X_t^{-1}X_L^{1/2}||_2 \,||X_L^{-1}||_F\\
     &\leq \left( -1-W_{-1}\left(-e^{-1-\epsilon_c}\right)\right)||X_L^{-1}||_F\\
     &\leq \left(\epsilon_c + \sqrt{2\epsilon_c}\right)||X_L^{-1}||_F
\end{split}
\label{app:eq:Tbounds}
\end{equation}
where the second inequality comes from the bounds derived in \cite{Roig2020} based on the Lambert $W\left(\cdot\right)$ function for maximum determinant problems, and the last inequality follows from the bound $W_{-1}\left(-e^{-1-u}\right)\geq -1-\sqrt{2u}-u$ derived in \cite{Chatz2013}. The bound on $t$ becomes:
\begin{equation}
     || E||_F +\left(\epsilon_c + \sqrt{2+\epsilon_c}\right)||X_L^{-1}||_F \geq |-\tau-t|
     \label{App:eq:Prebound}
\end{equation}
The norm on the residual $E$ can be bounded in the following way. Instantiate the optimization problem \eqref{4:eq:PSDcent} using the Cholesky bases $U_L$ of $X_L$ and compute a Newton step $\Delta_Y$ from the feasible point $Y=I$. This leads to:
\begin{equation}
    -\nabla^{2}h(I) \Delta_{Y} = \nabla h(I) + \tilde{\tau} U_L C U_L^T \, + \,\sum_i^M \tilde{\beta_i} U_L A_i U_L^T 
\end{equation}
where $\tilde{\tau}$ and $\tilde{\beta_i}$ are the dual variables corresponding to the linear equality constraints of \eqref{4:eq:PSDcent}. Evaluating the gradient and Hessian of $h(I)$ leads to:
\begin{equation}
     \Delta_{Y} = I - \frac{\tilde{\tau}}{N-1} U_L C U_L^T \, - \,\sum_i^M \frac{\tilde{\beta_i}}{N-1} U_L A_i U_L^T 
\end{equation}
And multiplying left and right by $U_L^{-1}$ and $U_L^{-T}$ respectively yields:
\begin{equation}
    U_L^{-1} \Delta_{Y} U_L^{-T}= X_L^{-1} - \frac{\tilde{\tau}}{N-1}  C  \, - \,\sum_i^M \frac{\tilde{\beta_i}}{N-1}  A_i 
\end{equation}
where we can draw the equivalence $E = U_L^{-1} \Delta_{Y} U_L^{-T}$, $\tau = \frac{\tilde{\tau}}{N-1}$ and $\beta_i =  \frac{\tilde{\beta_i}}{N-1}$ comparing the last equation to \eqref{App:eq:tboundKKT}. The norm of $E$ can then be bounded as:
\begin{equation}
    ||E||_F = ||U_L^{-1} \Delta_{Y} U_L^{-T}||_F \leq || \Delta_{Y} ||_F||U_L^{-1} U_L^{-T}||_F \leq \lambda_h(I) ||X_L^{-1}||F
\end{equation}
where $\lambda_h(\cdot)$ is the Newton decrement of problem \eqref{4:eq:PSDcent} and we use an inequality on Frobenius norms and the properties of the Newton step for self-concordant functions. Using the self-concordance of $-h(Y)$ and the optimality bound $\epsilon_C$, we have that:
\begin{equation}
    -h(I) -\alpha\beta \lambda_h(I) \geq -h(Y^*) \geq -h(I) - \epsilon_C
\end{equation}
where $\alpha\in(0,0.5)$ and $\beta\in(0,1)$ are the line-search parameters used in Section \ref{Subsec:Centering} and leads to $\frac{\epsilon_C}{\alpha\beta}\geq \lambda_h(I)$ and finally to:
\begin{equation}
    ||E||_F \leq \frac{\epsilon_C}{\alpha\beta} ||X_L^{-1}||F
\end{equation}
Combining this result with \eqref{App:eq:Prebound} leads to:
\begin{equation}
     \left(\left(1+\frac{1}{\alpha\beta}\right)\epsilon_c + \sqrt{2+\epsilon_c}\right)||X_L^{-1}||_F \geq |-\tau-t|
\end{equation}
leading to the final bounds $t^+ \geq t \geq t^-$:
\begin{equation}
    \begin{split}
        t^+ &= -\tau +  \left(\left(1+\frac{1}{\alpha\beta}\right)\epsilon_c + \sqrt{2\epsilon_c}\right)||X_L^{-1}||_F\\
        t^- &= -\tau -   \left(\left(1+\frac{1}{\alpha\beta}\right)\epsilon_c + \sqrt{2\epsilon_c}\right)||X_L^{-1}||_F
    \end{split}
\end{equation}
\end{proof}
\section{Gradients and Hessians of the Conic Barrier Functions}\label{Appendix:GradsHessians}
 The gradients of $\phi_{DD}$ and $\phi_{SDD}$ are separable, with each block of the  gradient with respect to $m_{i,j} = [x,y,z]$  given by:
\begin{equation}
    \nabla_{\phi_{DD}}(m_{i,j}) = 
\begin{bmatrix}
    \frac{x}{x^2-z^2} \\
    \frac{y}{y^2-z^2} \\
    z\left(\frac{-1}{x^2-z^2} + \frac{-1}{y^2-z^2}  \right)
    \end{bmatrix}; \;
    \nabla_{\phi_{SDD}}(m_{i,j}) = \frac{1}{xy-z^2}
\begin{bmatrix}
    y\\
    x\\
    -2z
    \end{bmatrix},
\end{equation}
Similarly, the $3 {N \choose 2} \times 3 {N \choose 2}$ Hessians are block-diagonal matrix composed of $N \choose 2$ blocks of size $3\times3$. Each block is indexed by a pair $\{i,j\}$ and is of the form:
\begin{equation} \begin{split}
&    \nabla^2_{\phi_{DD}}(m_{i,j}) = 
 \begin{bmatrix}
    \frac{-x^2-z^2}{(x^2-z^2)^2} & 0 &\frac{2xz}{(x^2-z^2)^2}\\
    *&\frac{-y^2-z^2}{(y^2-z^2)^2} &\frac{2yz}{(y^2-z^2)^2} \\
     * & *& \frac{-x^2-z^2}{(x^2-z^2)^2}+\frac{-y^2-z^2}{(y^2-z^2)^2}
    \end{bmatrix} \\
   & \nabla^2_{\phi_{SDD}}(m_{i,j}) = 
\begin{bmatrix}
   \frac{-y^2}{(xy-z^2)^2} & \frac{-z^2}{(xy-z^2)^2} & \frac{2yz}{(xy-z^2)^2}\\
    *&\frac{-x^2}{(xy-z^2)^2} & \frac{2xz}{(xy-z^2)^2} \\
     * & *& \frac{-2(xy+z^2)}{(xy-z^2)^2}
    \end{bmatrix}
    \end{split}
\end{equation}
\noindent The gradient of  $h(Y)$
is given by $\nabla_h(Y) = \left(N-1\right)\,Y^{-1}$. Its Hessian is a 4-dimensional tensor, whose quadratic form can be expressed as $-\left(N-1\right)\,\operatorname{Tr}(Y^{-1}UY^{-1}V)$, where $U$ and $V$ are the arguments of the quadratic form \cite{Boyd2004}.

\end{document}